\documentclass[11pt,a4paper,reqno]{amsart}

\usepackage[utf8x]{inputenc}

\usepackage{stmaryrd}
\usepackage{amssymb}
\usepackage{amsmath}
\usepackage{amsfonts}
\usepackage{amsthm}
\usepackage{color}
\usepackage[usenames,dvipsnames]{xcolor} 
\usepackage{bbm}
\usepackage{graphicx}
\usepackage{rotating}
\usepackage{hyperref}
\usepackage{mathrsfs}
\usepackage{cmll}
\usepackage[all, cmtip]{xy}
\usepackage{amscd}
\usepackage{tikz-cd}
\usepackage{todonotes}
\usepackage{tensor}
\providecommand{\U}[1]{\protect\rule{.1in}{.1in}}
\newtheorem{theorem}{Theorem}[subsection]
\newtheorem{corollary}[theorem]{Corollary}
\newtheorem{lemma}[theorem]{Lemma}
\newtheorem{proposition}[theorem]{Proposition}
\newtheorem*{theorem*}{Theorem}

\theoremstyle{definition}
\newtheorem{definition}[theorem]{Definition}
\newtheorem{examplex}[theorem]{Example}

\theoremstyle{remark}
\newtheorem{remarkx}[theorem]{Remark}
\newenvironment{remark}
  {\pushQED{\qed}\remarkx}
  {\popQED\endremarkx}

\thanks{}

\numberwithin{equation}{section}

\makeatletter
\def\@tocline#1#2#3#4#5#6#7{\relax
  \ifnum #1>\c@tocdepth 
  \else
    \par \addpenalty\@secpenalty\addvspace{#2}%
    \begingroup \hyphenpenalty\@M
    \@ifempty{#4}{%
      \@tempdima\csname r@tocindent\number#1\endcsname\relax
    }{%
      \@tempdima#4\relax
    }%
    \parindent\z@ \leftskip#3\relax \advance\leftskip\@tempdima\relax
    \rightskip\@pnumwidth plus4em \parfillskip-\@pnumwidth
    #5\leavevmode\hskip-\@tempdima
      \ifcase #1
       \or\or \hskip 1.5 em \or \hskip 2em \else \hskip 3em \fi%
      #6\nobreak\relax
    \hfill\hbox to\@pnumwidth{\@tocpagenum{#7}}\par
    \nobreak
    \endgroup
  \fi}
\makeatother

\begin{document}

\newcommand{\R}{{\mathbbm R}}
\newcommand{\C}{{\mathbbm C}} 
\newcommand{\T}{{\mathbbm T}}
\newcommand{\D}{{\mathbbm D}}
\renewcommand{\P}{\mathbb P}

\newcommand{\Aa}{{\mathcal A}}
\newcommand{\Ii}{{\mathbbm K}}
\newcommand{\Jj}{{\mathbbm J}}
\newcommand{\Nn}{{\mathcal N}}
\newcommand{\Ll}{{\mathcal L}}
\newcommand{\Tt}{{\mathcal T}}
\newcommand{\Gg}{{\mathcal G}}
\newcommand{\Dd}{{\mathcal D}}
\newcommand{\Cc}{{\mathcal C}}
\newcommand{\Oo}{{\mathcal O}}

\newcommand{\Hom}            {\operatorname{\mathsf{Hom}}}

\newcommand{\pr}{\operatorname{pr}}
\newcommand{\bla}{\langle \hspace{-2.7pt} \langle}
\newcommand{\bra}{\rangle\hspace{-2.7pt} \rangle}
\newcommand{\blq}{[ \! [}
\newcommand{\brq}{] \! ]}
 \newcommand{\into}{\mathbin{\vrule width1.5ex height.4pt\vrule height1.5ex}}

\newcommand{\TODO}[1]{
     \todo[
     size=\tiny]{#1}
     }

\setlength{\marginparwidth}{2cm}

\title{Deformations of vector bundles over Lie groupoids}

\author{Pier Paolo La Pastina}
\address{Dipartimento di Matematica ``Guido Castelnuovo'', Universit\`a degli studi di Roma ``La Sapienza'', P.le Aldo Moro 5, I-00185 Roma, Italy.}
\email{lapastina@mat.uniroma1.it}

\author{Luca Vitagliano}
\address{DipMat, Universit\`a degli Studi di Salerno, via Giovanni Paolo II n${}^\circ$ 123, 84084 Fisciano (SA) Italy.}
\email{lvitagliano@unisa.it}

\begin{abstract}
VB-groupoids are vector bundles in the category of Lie groupoids. They encompass several classical objects, including Lie group representations and 2-vector spaces. Moreover, they provide geometric pictures for 2-term representations up to homotopy of Lie groupoids. We attach to every VB-groupoid a cochain complex controlling its deformations and discuss its fundamental features, such as Morita invariance and a van Est theorem. Several examples and applications are given.
\end{abstract}

\maketitle

\tableofcontents

\section*{Introduction}

Lie groupoids are unifying structures in differential geometry. Their theory generalizes that of Lie groups: they have an infinitesimal counterpart, Lie algebroids, and there exists a Lie differentiation functor from Lie groupoids to Lie algebroids. Moreover, Lie groupoids provide a general framework to deal with many geometric situations such as orbifolds, foliations, Lie group actions and the geometry of PDEs. In all these cases, the main difficulty lies in studying spaces that are obtained as quotients of smooth manifolds, but are not smooth: they are rather \emph{differentiable stacks} \cite{behrend:diff} and it is well-known that the latter are morally just Lie groupoids up to Morita equivalence.

This leads to the idea that the category of Lie groupoids should be understood as a more general setting for differential geometry, and indeed several kinds of geometric structures on Lie groupoids (symplectic structures, complex structures, Riemannian metrics and so on) have been considered in the last years. In particular, vector bundles in the category of Lie groupoids are known as \emph{VB-groupoids}. They first appeared in \cite{pradines:theorie} in connection to symplectic groupoids. Later, in \cite{gracia-saz:vb2}, it was shown that they are geometric and more intrinsic models for 2-term representations up to homotopy of Lie groupoids \cite{abad:ruth2}. For example, the tangent and the cotangent bundle of a Lie groupoid are VB-groupoids: they correspond to the adjoint and the coadjoint representations.

This paper is part of a bigger project devoted to the study of deformations of geometric structures on Lie groupoids (and on differentiable stacks). Many such structures can be understood as vector bundle maps, so it is important to understand deformations of VB-groupoids first: we address this problem in the present paper concentrating on cohomological aspects. Our starting point is the paper \cite{crainic:def2} by Crainic, Mestre and Struchiner, where deformations of Lie groupoids are studied. We will also in part follow \cite{lapastina:def} where we discussed deformations of VB-algebroids, the infinitesimal version of VB-groupoids.

The paper is divided in two main sections: the first one presents the relevant deformation cohomology and its main properties, the second one discusses examples and applications. In its turn, the first section is divided in six subsections. In Subsection \ref{Sec:deformations1} we recall from \cite{crainic:def2, crainic:def, lapastina:def} the necessary results about deformations of Lie groupoids, Lie algebroids and VB-algebroids, while in Subsection \ref{Sec:deformations2} we recall the basic facts about VB-groupoids. 
In Subsection \ref{sec:lin_def} we discuss the deformation theory of VB-groupoids. If $(\mathcal V \rightrightarrows E; \mathcal G \rightrightarrows M)$ is a VB-groupoid, the top groupoid $\mathcal V \rightrightarrows E$ has an associated deformation complex $C_{\mathrm{def}}(\mathcal V)$. We show that (infinitesimal) deformations of the VB-groupoid structure are controlled by a subcomplex $C_{\mathrm{def,lin}}(\mathcal V)$ of $C_{\mathrm{def}}(\mathcal V)$, the \emph{linear deformation complex} of $\mathcal V$, originally introduced in \cite{etv:infinitesimal}. We compute the first cohomology groups and we show that dual VB-groupoids have the same linear deformation cohomology. In this section we will need a technical result, that is proved in Appendix \ref{sec:app}. In Subsection \ref{sec:linearization_1} we introduce the \emph{linearization map}, an important technical tool, adapting ideas from \cite{cabrera:hom} and \cite{lapastina:def}. The linearization map is used to prove that the linear deformation cohomology of a VB-groupoid is embedded in the deformation cohomology of the top groupoid. The linearization map also has applications in the subsequent subsections.

Subsection \ref{sec:van_est} is dedicated to the proof of a van Est theorem for the linear deformation cohomology of VB-groupoids. We show that the linear deformation complex of a VB-groupoid and that of the associated VB-algebroid are intertwined by a van Est map, which is a quasi-isomorphism under certain connectedness conditions. Finally, in Subsection \ref{sec:morita} we show that the linear deformation cohomology is Morita invariant. This is particularly important, because it means that this cohomology is really an invariant of the associated vector bundle of differentiable stacks.

In the second section we deal with examples. Subsection \ref{sec:VB_grp} is about vector bundles in the category of Lie groups. The latter are equivalent to Lie group representations and  we show that, in this case, the linear deformation cohomology and the classical cohomology controlling deformations of Lie group representations  \cite{nijenhuis:def} fit into a long exact sequence. In Subsection \ref{sec:2-vect} we study deformations of 2-vector spaces, i.e. Lie groupoids in the category of vector spaces. In Subsection \ref{sec:tangent-VB}, we discuss deformations of the tangent and the cotangent VB-groupoids of a Lie groupoid. Representations of foliation groupoids and Lie group actions on vector bundles can be also encoded by VB-groupoids: we study the associated deformation complexes in Subsections \ref{sec:fol_grpd} and \ref{sec:Lie_vect} respectively.

We assume that the reader is familiar with the basic notions on Lie groupoids and algebroids. Despite they are not as standard objects as the latter, we will also assume familiarity with VB-algebroids. The unfamiliar reader may refer e.g.~to \cite{gracia-saz:vb, bursztyn:vec} for the details, and \cite{lapastina:def} for our notation/conventions about VB-algebroids.

Before we start, we recall from \cite{grab:vect} the concept of \emph{homogeneity structure} of a vector bundle, a basic tool that will be used throughout the paper. Let $E \to M$ be a vector bundle. The monoid $\mathbb R_{\geq 0}$ of non-negative real numbers $\lambda$ acts on $E$ by homotheties $h_\lambda: E \to E$ (fiber-wise scalar multiplication). The action $h : \mathbb R_{\geq 0} \times E \to E$, $e \mapsto h_\lambda (e)$, is called the \emph{homogeneity structure} of $E$. Together with the smooth structure, it fully characterizes the vector bundle structure. This implies that every notion that involves the linear structure of $E$ can be expressed in terms of $h$ only: for example, a smooth map between the total spaces of two vector bundles is a vector bundle map if and only if it commutes with the homogeneity structures.

\section{Deformations of VB-groupoids}

\subsection{Preliminaries}\label{Sec:deformations1}

In this section, we briefly recall the deformation theory of Lie groupoids, Lie algebroids and VB-algebroids, originally presented in \cite{crainic:def2}, \cite{crainic:def} and \cite{lapastina:def} respectively.

\subsubsection{Deformations of Lie groupoids}

Let $\mathcal G \rightrightarrows M$ be a Lie groupoid. We denote by $\mathsf s, \mathsf t, \mathsf 1, \mathsf m, \mathsf i$ its structure maps (source, target, unit, multiplication, inversion respectively), by $\bar{\mathsf m}$ the division map and by $\mathcal{G}^{(k)}$ the manifold of $k$-tuples of composable arrows of $\mathcal{G}$. 

\begin{definition}\label{def:def_complex}
The \emph{deformation complex} $(C_{\mathrm{def}}(\mathcal{G}), \delta)$ of $\mathcal{G}$ is defined as follows. For $k \geq 0$, $C^k_{\mathrm{def}}(\mathcal{G})$ is the set of smooth maps $c: \mathcal{G}^{(k+1)} \to T\mathcal{G}$ such that:
\begin{enumerate}
	\item $c(g_0, \dots, g_k) \in T_{g_0} \mathcal{G}$;
	\item $(T \mathsf s \circ c) (g_0, \dots, g_k)$ does not depend on $g_0$
\end{enumerate}
for each $(g_0, \dots, g_k) \in \mathcal{G}^{(k+1)}$. Thus we define the \emph{$\mathsf s$-projection of $c$} to be 
\[
s_c: \mathcal G^{(k)} \to TM, \quad s_c(g_1, \dots, g_k) := (T\mathsf s \circ c)(g_0, \dots, g_k).
\]

The differential of $c \in C^k_{\mathrm{def}} (\mathcal{G})$ is defined by
\begin{equation}\label{eq:diff}
\begin{aligned}
\delta c (g_0, \dots, g_{k+1}) = & -T \bar{\mathsf m} (c(g_0 g_1, \dots, g_{k+1}), c(g_1, \dots, g_{k+1})) \\ 
& + \sum_{i=1}^k (-1)^{i-1} c(g_0, \dots, g_i g_{i+1}, \dots, g_{k+1}) + (-1)^k c(g_0, \dots, g_k).
\end{aligned}
\end{equation}
Moreover, $C^{-1}_{\mathrm{def}}(\mathcal{G}) := \Gamma(A)$, where $A \Rightarrow M$ is the Lie algebroid of $\mathcal{G}$, and $\delta \alpha = \overleftarrow{\alpha} + \overrightarrow{\alpha}$ for each $\alpha \in \Gamma(A)$, where $\overleftarrow{\alpha}$ and $\overrightarrow{\alpha}$ are the left-invariant and right-invariant vector fields determined by $\alpha$.
\end{definition}

\begin{remark} Notice that we adopt a different convention from \cite{crainic:def2}, where $C^k_{\mathrm{def}}(\mathcal G)$ is the space of smooth maps $\mathcal G^{(k)} \to T\mathcal G$ satisfying properties (1) and (2), to be coherent with \cite{lapastina:def}, where deformations of VB-algebroids were studied.  \end{remark}

\begin{remark}\label{rmk:couple} We observe, for later use, that conditions (1) and (2) can be expressed in the following way. For every $k \geq 1$, define the surjective submersions
\begin{equation}\label{eq:pq}
\begin{aligned}
p_k: \mathcal G^{(k)} \to \mathcal G, \quad & (g_1, \dots, g_k) \mapsto g_1, \\
q_k: \mathcal G^{(k)} \to M, \quad & (g_1, \dots, g_k) \mapsto \mathsf s(g_1).
\end{aligned}
\end{equation}
Then an element $c \in C^k_{\mathrm{def}}(\mathcal G)$ is simply a section of the pull-back bundle $p^*_{k+1} T \mathcal G \to \mathcal G^{(k+1)}$ that is $\mathsf s$-projectable, i.e. such that there exists a section $s_c$ of $q_k^* TM \to \mathcal G^{(k)}$ fitting into the following commutative diagram
\[
\begin{array}{c}
\xymatrix{p_{k+1}^* T \mathcal G \ar[d]_{T\mathsf s} \ar[r]& \mathcal G^{(k+1)} \ar@/_1.2pc/[l]_{c} \ar[d] &  \\
q_k^* TM \ar[r] & \mathcal G^{(k)} \ar@/_1.2pc/[l]_{s_c} }
\end{array}
,
\]
where the right arrow is the projection onto the last $k$ elements.
\end{remark}

Recall that there is a short exact sequence of vector bundles over $\mathcal G$:
\begin{equation}\label{eq:ses1}
0 \longrightarrow T^{\mathsf s} \mathcal G \longrightarrow T \mathcal G \overset{T\mathsf s}\longrightarrow {\mathsf s}^* TM \longrightarrow 0.
\end{equation}
Pulling back the sequence via $\mathsf 1: M \to \mathcal G$, we get a short exact sequence of vector bundles over $M$
\begin{equation}\label{eq:ses2}
0 \longrightarrow A \longrightarrow \mathsf 1^* T \mathcal G \longrightarrow TM \to 0.
\end{equation}
The latter has a canonical splitting, given by $T\mathsf 1: TM \to T \mathcal G$, so $\mathsf 1^* T \mathcal G \cong A \oplus TM$ canonically.

\begin{definition}\label{def:norm_complex}
The \emph{normalized deformation complex} $\widehat C_{\mathrm{def}}(\mathcal G)$ of $C_{\mathrm{def}}(\mathcal G)$ is defined as follows. For $k \geq 1$, $\widehat C^k_{\mathrm{def}}(\mathcal G)$ is composed by all cochains $c \in C^k_{\mathrm{def}}(\mathcal G)$ such that
\begin{equation}\label{eq:norm}
c(\mathsf 1_x, g_1, \dots, g_k) \in T_x M \subset T_{\mathsf 1_x} \mathcal G, \quad \text{and} \quad c(g_0, \dots, \mathsf 1_x, \dots, g_k) = 0,
\end{equation}
while $\widehat C^0_{\mathrm{def}}(\mathcal G)$ is composed by $0$-cochains $c$ that satisfy $c(\mathsf 1_x) = s_c(x)$. Finally, $\widehat C^{-1}_{\mathrm{def}}(\mathcal G) := \Gamma(A)$. 
\end{definition}

Notice that the first condition in (\ref{eq:norm}) implies that $c(\mathsf 1_x, g_1, \dots, g_k)$ can be identified with $(T\mathsf s \circ c)(\mathsf 1_x, g_1, \dots, g_k) = s_c (g_1, \dots, g_k)$, so we recover the definition in \cite{crainic:def2}. Moreover, we have:

\begin{proposition}[{\cite[Proposition 11.8]{crainic:def2}}]\label{prop:norm} The inclusion $\widehat C_{\mathrm{def}}(\mathcal G) \hookrightarrow C_{\mathrm{def}}(\mathcal G)$ is a quasi-isomorphism. \end{proposition}

The group of automorphisms of $\mathcal G$ acts naturally on $C_{\mathrm{def}}(\mathcal G)$ by pullback. Explicitly, if $\Psi \in \mathrm{Aut}(\mathcal G)$, denote by $\psi$ the induced automorphism of $A$: this notation will be used throughout the paper. Then the action is given by
\[
(\Psi^* c)(g_0, \dots, g_k) := T \Psi^{-1} (c(\Psi(g_0), \dots, \Psi(g_k)))
\]
for $c \in C^k_{\mathrm{def}}(\mathcal G), k \geq 0$, and by
\[
\Psi^* c := \psi^* c
\]
for $c \in C^{-1}_{\mathrm{def}}(\mathcal G)$. It is easy to check that this action preserves the differential. Indeed, if $\alpha \in \Gamma(A)$, we have
\begin{equation}\label{eq:inv_vf}
\Psi^* \overrightarrow{\alpha} = \overrightarrow{\psi^* \alpha}, \quad \Psi^* \overleftarrow{\alpha} = \overleftarrow{\psi^* \alpha}
\end{equation}
and so
\[
\Psi^* (\delta \alpha) = \Psi^* (\overleftarrow{\alpha} + \overrightarrow{\alpha}) = \overleftarrow{\psi^* \alpha} + \overrightarrow{\psi^* \alpha} = \delta (\psi^* \alpha) = \delta (\Psi^* \alpha).
\]
If $c \in C^k_{\mathrm{def}}(\mathcal G)$, $k \geq 0$, a direct computation exploiting Equation (\ref{eq:diff}) shows that $\Psi^*$ commutes with $\delta$. Finally, it is straightforward that the action preserves $\widehat C_{\mathrm{def}}(\mathcal G)$.

Recall that a \emph{representation} of $\mathcal{G}$ on a vector bundle $E \to M$ is a morphism from $\mathcal{G}$ to the general linear groupoid $GL(E)$. The isomorphism $E_x \to E_y$ induced by an arrow $g: x \to y$ of $\mathcal{G}$ is denoted $v \mapsto g \cdot v$.

\begin{definition} \label{def:coh-rep}
The \emph{(differentiable) complex of $\mathcal{G}$ with coefficients in the representation $E$}, denoted $C (\mathcal{G},E)$, is defined as follows. For $k > 0$, $C^k(\mathcal{G},E)$ is the set of smooth maps $u: \mathcal{G}^{(k)} \rightarrow E$ such that $u(g_1, \dots, g_k) \in E_{\mathsf t(g_1)}$. The differential of $u \in C^k(\mathcal{G},E)$ is given by
\[
\begin{aligned}
\delta u (g_1, \dots, g_{k+1}) =  & g_1 \cdot u(g_2, \dots, g_{k+1}) + \\ 
& \sum_{i=1}^k (-1)^i u(g_1, \dots, g_i g_{i+1}, \dots, g_{k+1}) + (-1)^{k+1} u(g_1, \dots, g_k).
\end{aligned}
\]
Moreover, $C^{0}(\mathcal G, E) := \Gamma(E)$, and, for $\varepsilon \in C^{0}(\mathcal G, E)$, $\delta \varepsilon$ is defined by
\[
\delta \varepsilon (g) = g \cdot \varepsilon_{\mathsf s(g)} - \varepsilon_{\mathsf t(g)}.
\]
The cohomology of $C (\mathcal G, E)$ is called the \emph{(differentiable) cohomology of $\mathcal{G}$ with coefficients in $E$} and denoted by $H (\mathcal G, E)$.

When $E$ is the trivial line bundle with the trivial representation, the above complex is simply called the \emph{Lie groupoid complex} of $\mathcal{G}$ and denoted $C (\mathcal{G})$. Its cohomology is called the \emph{Lie groupoid cohomology} of $\mathcal G$ and denoted by $H (\mathcal G)$. The complex $C (\mathcal G)$ possesses a canonical DG-algebra structure. The product of $f_1 \in C^k (\mathcal G)$ and $f_2 \in C^l(\mathcal G)$ is defined by
\begin{equation}\label{eq:prod}
(f_1 \star f_2)(g_1, \dots, g_{k+l}) = f_1 (g_1, \dots, g_k) f_2 (g_{k+1}, \dots, g_{k+l})
\end{equation}
if $k, l > 0$, by
\[
(f_1 \star f_2)(g_1, \dots, g_l) = f_1 (\mathsf t(g_1)) f_2 (g_1, \dots, g_l)
\]
if $k = 0$, $l > 0$, by
\[
(f_1 \star f_2)(g_1, \dots, g_k) = f_1 (g_1, \dots, g_k) f_2(\mathsf s(g_k))
\]
if $k > 0, l = 0$, and by
\[
f_1 \star f_2 = f_1 f_2
\]
if $k = l = 0$. The same formulas define a $C (\mathcal G)$-DG-module structure on $C_{\mathrm{def}}(\mathcal G)$ and on $C (\mathcal G, E)$, for every representation $E$.
\end{definition}

Now we recall from \cite{crainic:def2} the deformation cohomology of a Lie groupoid in degrees $-1$ and $0$. To do this, we need to recall the canonical representations on the isotropy and the normal bundles first. When the groupoid is regular, i.e.~its orbits have all the same dimension, the latter are smooth representations. In general, they are only set-theoretic representations, but we can still define invariant sections.

Let $\mathcal{G} \rightrightarrows M$ be a Lie groupoid, $A \Rightarrow M$ its Lie algebroid. The \emph{isotropy bundle} of $\mathcal{G}$ is defined by
\[
\mathfrak{i} := \ker (\rho: A \to TM).
\]
The bracket of $A$ induces a Lie algebra structure on each fiber $\mathfrak i_x$ of $\mathfrak{i}$, $x \in M$. The Lie algebra $\mathfrak i_x$ is actually the Lie algebra of the isotropy group $\mathcal G_x$ of $\mathcal G$ at $x$. For every $g: x \to y$ in $\mathcal{G}$, there is an obvious conjugation map $Ad_g: \mathcal{G}_x \to \mathcal{G}_y$. Differentiating at units one obtains the \emph{isotropy action} of $\mathcal{G}$ on $\mathfrak{i}$, 
\[
ad_g: \mathfrak{i}_x \to \mathfrak{i}_y.
\]
Even if $\mathfrak{i}$ is not a honest vector bundle, we define sections of $\mathfrak{i}$ to be
\[
\Gamma(\mathfrak{i}) := \ker (\rho: \Gamma(A) \to \mathfrak{X}(M)),
\]
and invariant sections
\[
H^0(\mathcal{G}, \mathfrak{i}) = \Gamma(\mathfrak{i})^{\mathrm{inv}} := \{ \alpha \in \Gamma(\mathfrak{i}) : ad_g (\alpha_x) = \alpha_y \ \forall \ g: x \to y \ in \ \mathcal{G} \}.
\]

\begin{proposition}\label{prop:H^-1} $H^{-1}_{\mathrm{def}}(\mathcal{G}) \cong H^0 (\mathcal{G}, \mathfrak{i}) = \Gamma(\mathfrak{i})^{\mathrm{inv}}$. \end{proposition}

Again, even if $\mathfrak{i}$ is not of constant rank, one can define the differentiable cohomology of $\mathcal{G}$ with coefficients in $\mathfrak{i}$ as in Definition \ref{def:coh-rep}, where a cochain is smooth if it is smooth as an $A$-valued map. One can prove that the complex obtained in this way is well defined and there is an inclusion of complexes
\begin{equation}\label{eq:isotropy}
r: C (\mathcal G, \mathfrak i) \hookrightarrow C_{\mathrm{def}} (\mathcal G)[-1] , \quad u \mapsto c_u
\end{equation}
given by
\[
c_u(g_1, \dots, g_k) = T R_{g_1} (u(g_1, \dots, g_k)).
\]
where $R_g$ denotes \emph{right translation} by $g \in \mathcal G$.

Now recall that, for any Lie groupoid $\mathcal{G} \rightrightarrows M$, one can construct the \emph{tangent prolongation Lie groupoid} $T\mathcal{G} \rightrightarrows TM$: its structure maps are the tangent maps to the original ones. A vector field $X \in \mathfrak X (\mathcal G)$ is called \emph{multiplicative} if it is a groupoid morphism $X: \mathcal G \to T \mathcal G$. Notice that, if $\alpha \in \Gamma (A)$, $\delta \alpha = \overrightarrow{\alpha} + \overleftarrow{\alpha}$ is a multiplicative vector field. The flow of $\delta \alpha$ consists of inner automorphisms of $\mathcal G$ and every $\delta \alpha$ is called an \emph{inner multiplicative vector field}.

\begin{proposition}\label{prop:H^0}
\[
H^0_{\mathrm{def}}(\mathcal G) = \dfrac{\text{multiplicative vector fields on $\mathcal G$}}{\text{inner multiplicative vector fields on $\mathcal G$}}
\]
\end{proposition}

The next step will be the definition of the normal representation. The \emph{normal bundle} of $\mathcal G$ is defined by
\[
\nu := TM/ \operatorname{im} \rho.
\]
Even if $\nu$ is not a smooth vector bundle, we can define its sections by
\begin{equation}\label{eq:sections_nu}
\Gamma(\nu) := \dfrac{\mathfrak X (M)}{\operatorname{im} (\rho: \Gamma(A) \to \mathfrak X (M))}.
\end{equation}
There is a natural action of $\mathcal G$ on $\nu$, defined as follows. Take an arrow $g: x \to y$ in $\mathcal G$ and $v \in \nu_x$, then choose a curve $g(\epsilon): x(\epsilon) \to y(\epsilon)$ such that $g(0) = g$ and $\dot x (0)$ represents $v$; we define $ad_g (v) = [\dot y (0)]$. One can check that this definition does not depend on the choices involved and that the following lemma holds.

\begin{lemma} Let $V \in \mathfrak X (M)$. Then the set-theoretic section of $\nu$
\[
x \mapsto [V_x]
\]
is invariant if and only if, for any $\mathsf s$-lift $X$ of $V$ and any $g: x \to y$, there exists $\eta (g) \in A_{\mathsf t(g)}$ such that
\begin{equation}\label{eq:invariance}
V_{\mathsf t(g)} = T\mathsf t(X_g) + \rho (\eta (g)).
\end{equation}
\end{lemma}
Moreover, if Equation (\ref{eq:invariance}) holds for some $\mathsf s$-lift, it holds for all of them.

Now notice that every section of $\nu$, as defined in (\ref{eq:sections_nu}), induces a set-theoretic section. So it is natural to declare that a section $[V]$ of $\nu$ is invariant if, for some $\mathsf s$-lift $X$ of $V$, there exists a smooth section of $\mathsf t^* A \to \mathcal G$ such that Equation (\ref{eq:invariance}) holds. But the vector field $X' \in \mathfrak X (\mathcal G)$ defined by
\[
X'_g = X_g + T R_g (\eta_g)
\]
is also a $\mathsf s$-lift, so we end up with the following definition.

\begin{definition} A section $[V] \in \Gamma(\nu)$ is \emph{invariant} if there exists $X \in \mathfrak X (\mathcal G)$ that is $\mathsf s$-projectable and $\mathsf t$-projectable to $V$. We say that $X$ is an \emph{$(\mathsf s, \mathsf t)$-lift} of $V$. The space of invariant sections is denoted $H^0 (\mathcal G, \nu)$ or $\Gamma(\nu)^{\mathrm{inv}}$. \end{definition}

There is a linear map
\begin{equation}\label{eq:pi}
\pi: H^0_{\mathrm{def}}(\mathcal G) \to \Gamma(\nu)^{\mathrm{inv}}
\end{equation}
that sends the class of a multiplicative vector field to the class of its projection on $M$.

\begin{lemma}[The curvature map] \label{prop:curvature} Let $[V] \in \Gamma(\nu)^{\mathrm{inv}}$ and $X$ an $(\mathsf s, \mathsf t)$-lift of $V$. Then $\delta X \in C^2(\mathcal G, \mathfrak i)$ and its cohomology class does not depend on the choice of $X$. Therefore there is an induced linear map
\begin{equation}\label{eq:curvature}
K: \Gamma(\nu)^{\mathrm{inv}} \to H^2(\mathcal G, \mathfrak i).
\end{equation}
\end{lemma}

Finally, we obtain:

\begin{proposition}\label{prop:exact} There is an exact sequence
\[
0 \longrightarrow H^1(\mathcal G, \mathfrak i) \overset{r}{\longrightarrow} H^0_{\mathrm{def}}(\mathcal G) \overset{\pi}{\longrightarrow} \Gamma(\nu)^{\mathrm{inv}} \overset{K}{\longrightarrow} H^2(\mathcal G, \mathfrak i) \overset{r}{\longrightarrow} H^1_{\mathrm{def}}(\mathcal G).
\]
\end{proposition}

The group $H^1_{\mathrm{def}}(\mathcal G)$ is directly linked to (infinitesimal) deformations. Before discussing this relation, we need a more general definition.

\begin{definition} A \emph{family of Lie groupoids}
\[
\tilde{\mathcal G} \rightrightarrows \tilde{M} \overset{\pi}{\longrightarrow} B,
\]
consists of a Lie groupoid $\tilde{\mathcal G} \rightrightarrows \tilde M$ and a surjective submersion $\pi: \tilde M \to B$ such that $\pi \circ \tilde {\mathsf s} = \pi \circ \tilde {\mathsf t}$. In particular, for every $b \in B$, $\mathcal G_b := (\pi \circ \tilde {\mathsf s})^{-1}(b)$ is a Lie groupoid over $M_b = \pi^{-1}(b)$. \end{definition}

This definition encodes the idea of a ``smoothly varying'' Lie groupoid. If $B$ is an open interval $I$ containing $0$, we say that $\tilde{\mathcal G}$ is a \emph{deformation} of $\mathcal G_0 \rightrightarrows M_0$ and we denote the latter by $\mathcal G \rightrightarrows M$. We will often denote by $\epsilon$ the canonical coordinate on $I$. Accordingly, a deformation of $\mathcal G$ is also denoted by $(\mathcal G_{\epsilon})$. The structure maps of $\mathcal G_{\epsilon}$ are denoted $\mathsf s_{\epsilon}, \mathsf t_{\epsilon}, \mathsf 1_{\epsilon}, \mathsf m_{\epsilon}, \mathsf i_{\epsilon}$. The division map is denoted $\bar {\mathsf m}_{\epsilon}$. A deformation $(\mathcal G_\epsilon)$ is called \emph{strict} if $\mathcal G_{\epsilon} \cong \mathcal G$ as manifolds for all $\epsilon$. This amounts to say that $\tilde{\mathcal G} \cong \mathcal G \times I$ and $\tilde M \cong M \times I$. In the following, for a strict deformation, we will always assume $\tilde{\mathcal G} = \mathcal G \times I$ (and $\tilde M = M \times I$). A strict deformation is \emph{$\mathsf s$-constant} (resp.~\emph{$\mathsf t$-constant}) if $\mathsf s_{\epsilon}$ (resp.~$\mathsf t_{\epsilon}$) does not depend on $\epsilon$. A deformation which is both $\mathsf s$-constant and $\mathsf t$-constant is \emph{$(\mathsf s, \mathsf t)$-constant}. The \emph{constant deformation} is the one with $\mathcal G_{\epsilon} = \mathcal G$ as groupoids for all $\epsilon$.

Two deformations $(\mathcal G_{\epsilon})$ and $(\mathcal G'_{\epsilon})$ of $\mathcal G$ are said to be \emph{equivalent} if there exists a smooth family of groupoid isomorphisms $\Psi_{\epsilon}: \mathcal G_{\epsilon} \to \mathcal G'_{\epsilon}$ such that $\Psi_0 = \mathrm{id}_{\mathcal G_0}$. We say that $(\mathcal G_{\epsilon})$ is \emph{trivial} if it is equivalent to the constant deformation.

Let $(\mathcal G_{\epsilon})$ be a strict deformation of the Lie groupoid $\mathcal G \rightrightarrows M$. Then it is natural to look at the variation of the structure maps, in particular the multiplication. However, in general, if $(g,h) \in \mathcal G^{(2)}$, there is no guaranty that $g$ and $h$ are composable also with respect to the groupoid structure $\mathcal G_{\epsilon}$. We will consider this problem later: now we simply assume that $(\mathcal G_{\epsilon})$ is an $(\mathsf s, \mathsf t)$-constant deformation. In this case, it makes sense to consider the tangent vector
\begin{equation}\label{eq:vector}
- \dfrac{d}{d \epsilon} \bigg|_{\epsilon = 0} \mathsf m_{\epsilon}(g,h) \in T_{gh} \mathcal G
\end{equation}
for any $(g,h) \in \mathcal G^{(2)}$. It is clear that (\ref{eq:vector}) is killed by both $T \mathsf s$ and $T \mathsf t$. This means that it is of the form $T R_{gh}(a)$ with $a \in A$, and moreover $a \in \mathrm{ker}(\rho) = \mathfrak i$. Hence we can define a cochain $u_0 \in C^1 (\mathcal G, \mathfrak i)$ by
\[
u_0 (g,h) := - \dfrac{d}{d \epsilon} \bigg|_{\epsilon = 0} R^{-1}_{gh} \left(\mathsf m_{\epsilon}(g,h)\right).
\]
Differentiating the associativity equation $\mathsf m_{\epsilon} (\mathsf m_{\epsilon}(g,h),k) = \mathsf m_{\epsilon}(g, \mathsf m_{\epsilon}(h,k))$ at $\epsilon = 0$, we find that $u_0$ is a cocycle.

Now, define $\xi_0 := r(u_0) \in C^1_{\mathrm{def}}(\mathcal G)$, where $r$ is the map (\ref{eq:isotropy}). Differentiating at 0 the identity $\mathsf m_{\epsilon}(\bar {\mathsf m}_{\epsilon} (\mathsf m_0(g,h),h),h) = \mathsf m_0(g,h)$, one obtains the following expression for $\xi_0$:
\begin{equation}\label{eq:def_coc}
\xi_0(g,h) = \dfrac{d}{d\epsilon} \bigg|_{\epsilon = 0} \bar {\mathsf m}_{\epsilon} (gh,h).
\end{equation}

The last computation suggests how to generalize the procedure. Let $(\mathcal G_{\epsilon})$ be an $\mathsf s$-constant deformation of $\mathcal G \rightrightarrows M$. The \emph{deformation cocycle} $\xi_0 \in C^1_{\mathrm{def}}(\mathcal G)$ associated to $(\mathcal G_{\epsilon})$ is defined by Formula (\ref{eq:def_coc}) (which makes sense also in the present case).

\begin{lemma}\label{prop:s-const} $\xi_0$ is a cocycle and its cohomology class only depends on the equivalence class of the deformation. \end{lemma}

Next we interpret $\xi_0$ in terms of the groupoid $\tilde{\mathcal G}$.

\begin{proposition} Let $\tilde{\mathcal G}$ be an $\mathsf s$-constant deformation of the Lie groupoid $\mathcal G$. Then, if we set $\xi = \delta (\frac{\partial}{\partial \epsilon}) \in C^1_{\mathrm{def}}(\tilde {\mathcal G})$, we have $\xi_0 = \xi|_{\mathcal G}$. \end{proposition}

Notice that this statement relies heavily on the fact that $(\mathcal G_{\epsilon})$ is $\mathsf s$-constant: otherwise the vector field $\frac{\partial}{\partial \epsilon}$ would not be $\tilde {\mathsf s}$-projectable. In the general case, one has to find an analogue of $\frac{\partial}{\partial \epsilon}$: this leads to the concept of \emph{transverse vector field}.

\begin{definition} Let $\tilde {\mathcal G}$ be a deformation of $\mathcal G$. A \emph{transverse vector field} for $\tilde{\mathcal G}$ is a vector field $X \in \mathfrak X (\tilde{\mathcal G})$ which is $\mathsf s$-projectable to a vector field $V \in \mathfrak X (\tilde M)$ which is, in turn, $\pi$-projectable to $\frac{d}{d \epsilon}$. \end{definition}

\begin{proposition}\label{prop:gen_def} Let $\tilde{\mathcal G}$ be a deformation of $\mathcal G$. Then:
\begin{enumerate}
\item there exist transverse vector fields for $\tilde{\mathcal G}$;
\item if $\tilde X$ is transverse, then $\delta \tilde{X}$, when restricted to $\mathcal G$, induces a cocycle $\xi_0 \in C^1_{\mathrm{def}}(\mathcal G)$;
\item the cohomology class of $\xi_0$ does not depend on the choice of $\tilde X$.
\end{enumerate}
\end{proposition}

The resulting cohomology class in $H^1_{\mathrm{def}}(\mathcal G)$ is called the \emph{deformation class} associated to the deformation $\tilde {\mathcal G}$. So, in general it is not possible to find a canonical cocycle. It was possible in the case of an $\mathsf s$-constant deformation because there was a canonical choice of a transverse vector field. Notice that from the proposition above it follows directly that the deformation class is also invariant under equivalence of deformations.

Finally, we recall a result about general families of Lie groupoids. Let $ \tilde{\mathcal G} \rightrightarrows \tilde M \overset{\pi}{\longrightarrow} B$ be a family of Lie groupoids. Then any curve $\gamma: I \to B$ induces a deformation $\gamma^* \tilde{\mathcal G}$ of $\tilde{\mathcal G}_{\gamma(0)}$. We have the following

\begin{proposition} [The variation map] Let $b \in B$. For any curve $\gamma: I \to B$ with $\gamma(0) = b$, the deformation class of $\gamma^* \tilde{\mathcal G}$ at time $0$ does only depend on $\dot{\gamma}(0)$. This defines a linear map
\[
\mathrm{Var}^{\tilde{\mathcal G}}_b: T_b B \to H^1_{\mathrm{def}}(\tilde{\mathcal G}_b),
\]
called the \emph{variation map} of $\tilde{\mathcal G}$ at $b$.
\end{proposition}

\subsubsection{Deformations of Lie algebroids and VB-algebroids}

Let $E \to M$ be a vector bundle. A \emph{multiderivation with $k$ entries} of $E$ (and $C^\infty (M)$-multilinear symbol), also called a $k$-\emph{derivation}, is a skew-symmetric, $\mathbb{R}$-$k$-linear map
\[
c: \Gamma(E) \times \dots \times \Gamma(E) \to \Gamma(E)
\]
such that there exists a bundle map $\sigma_c: \wedge^{k-1} E \to TM$, the \emph{symbol} of $c$, satisfying the following Leibniz rule:
\[
c(\varepsilon_1, \dots, \varepsilon_{k-1}, f \varepsilon_k) = \sigma_c (\varepsilon_1, \dots, \varepsilon_{k-1})(f) \varepsilon_k + f c(\varepsilon_1, \dots, \varepsilon_k),
\]
for all $\varepsilon_1, \dots, \varepsilon_k \in \Gamma(E)$, $f \in C^\infty (M)$. The space of $k$-derivations of $E \to M$ is denoted $\mathfrak D^k(E)$ and we set $\mathfrak D^\bullet (E) := \bigoplus_k \mathfrak D^k (E)$. 

$1$-derivations are simply derivations and they are of a particular interest. The space of derivations is denoted by $\mathfrak D (E)$. Recall that derivations of $E$ are sections of a Lie algebroid $DE \Rightarrow M$, that sits in the following short exact sequence (the \emph{Spencer sequence}):
\begin{equation} \label{eq:spencer}
0 \longrightarrow \operatorname{End} E \longrightarrow DE \overset{\sigma}{\longrightarrow} TM \longrightarrow 0,
\end{equation}
where $\sigma: DE \to TM$ is the \emph{symbol map} and plays the role of the anchor. Actually, $DE$ is the Lie algebroid of the general linear groupoid $GL(E)$ and it is called the \emph{gauge algebroid} of $E$.  Before going on, we recall its main properties. First of all, there is a canonical isomorphism of Lie algebroids 
\[
DE \overset{\cong}{\longrightarrow} DE^*, \quad \delta \mapsto \delta^*,
\]
defined as follows: if $\delta \in D_x E$, then $\delta^*: \Gamma(E^*) \to E^*_x$ is uniquely determined by the condition
\[
\sigma_\delta \langle \varphi, \varepsilon \rangle = \langle \delta^* \varphi, \varepsilon \rangle + \langle \varphi, \delta \varepsilon \rangle,
\]
for all $\varepsilon \in \Gamma(E)$, $\varphi \in \Gamma(E^*)$. 

There is an alternative description of derivations of a vector bundle $\pi : E \to M$. Namely, consider the vector bundle $T \pi: TE \to TM$. We denote by $TE|_v$ the fiber $T \pi^{-1} (v)$ over a tangent vector $v \in TM$ . Then a derivation $\delta \in D_x E$ (at a point $x \in M$) with symbol $\sigma_\delta$ determines a linear map $\widehat \delta: E_x \to TE|_{\sigma_{\delta}}$ via
\[
\widehat \delta (e) (\varphi) = \langle \delta^* \varphi, e \rangle,
\]
for all $e \in E_x$ and $\varphi \in \Gamma (E^\ast)$, where, in the lhs, $\varphi$ is also interpreted as a fiber-wise linear function on $E$. The assignment $\delta \mapsto (\sigma_\delta, \widehat \delta)$ establishes a one-to-one correspondence between the fiber $D_x E$ of the gauge algebroid over $x$, and the space of pairs $(v, h)$ where $v \in T_x M$ and $h : E_x \to TE|_v$ is a right inverse of the projection $TE|_v \to E_x$. In other words, $D E$ is the \emph{fat algebroid} \cite{gracia-saz:vb} of the VB-algebroid $(TE \Rightarrow E; TM \Rightarrow M)$, and derivations are the same as linear sections of $(TE \Rightarrow E; TM \Rightarrow M)$, or, which is the same, linear vector fields $\mathfrak X_{\mathrm{lin}} (E)$ on $E$: $\mathfrak D (E) \cong \mathfrak X _{\mathrm{lin}} (E)$ \cite{mackenzie}. We will sometimes identify $\delta$ and the pair $(\sigma_\delta, \widehat \delta)$.

The assignment $E \mapsto DE$ is functorial in the following sense. Let
\[
\begin{array}{c}
\xymatrix{E_N \ar[r]^{\phi} \ar[d] & E \ar[d] \\
N \ar[r]^{f} & M}
\end{array}
,
\]
be a \emph{regular} vector bundle morphism, i.e.~$\phi$ is an isomorphism on each fiber, so that $E_N$ is canonically isomorphic to the pullback bundle $ f^*E$. Then there is a pullback map $\phi^*: \Gamma(E) \to \Gamma(E_N)$ defined by
\[
(\phi^*\varepsilon)_y = \phi_y^{-1}(\varepsilon_{f (y)})
\]
for all $\varepsilon \in \Gamma (E)$, and all $y \in N$. One can use this pull-back to define a Lie algebroid morphism $D \phi: DE \to DF$. Specifically, for all $\delta \in D_y E_N$ we define $D \phi (\delta): \Gamma(E) \to E_{f(y)}$ by
\[
D \phi (\delta) = \phi \circ \delta \circ \phi^*.
\]
It is then easy to see that 
\begin{equation} \label{eq: delta_hat}
\widehat{D \phi (\delta)} = T \phi \circ \widehat \delta \circ \phi_x^{-1} : E_x \to TE|_{Tf (\sigma_{\delta})}.
\end{equation}
Finally, the diagram
\[
\xymatrix{
DE_N \ar[r]^-{D \phi} \ar[d] & DE \ar[d]\\
TN \ar[r]^-{Tf} & TM
}
\]
is a pull-back diagram, and this induces an isomorphism $DE_N \cong TN \times_{TM} DE$ which is sometimes useful. From now on, we will often identify $E_N$ with the pull-back $f^\ast E$.
For more details about the gauge algebroid we refer to \cite{mackenzie} (see also \cite{etv:infinitesimal}).

Now, let's go back to the main topic of this subsection and take a Lie algebroid $A \Rightarrow M$. We turn $C_{\mathrm{def}}(A) := \mathfrak D^\bullet (A)[1]$ into a cochain complex, the \emph{deformation complex} of $A$, using the Lie bracket $[-,-]$ on sections of $A \to M$. The differential $\delta: C^k_{\mathrm{def}} (A) \to C^{k+1}_{\mathrm{def}} (A)$ is defined by
\[
\begin{aligned}
\delta c (\alpha_0, \dots, \alpha_{k+1}) = & \sum_{i} (-1)^{i} [\alpha_i, c(\alpha_0, \dots, \hat{\alpha_i}, \dots, \alpha_{k+1})] \\
& + \sum_{i<j} (-1)^{i+j} c([\alpha_i, \alpha_j], \alpha_0, \dots, \hat{\alpha_i}, \dots, \hat{\alpha_j}, \dots, \alpha_{k+1}).
\end{aligned}
\]
for all $\alpha_0, \dots, \alpha_{k+1} \in \Gamma(A)$.

The group of Lie algebroid automorphisms of $A$ naturally acts on $C_{\mathrm{def}}(A)$ by pullback. Explicitly, if $\psi \in \mathrm{Aut}(A)$ and $c \in C^k_{\mathrm{def}}(A)$, then
\[
(\psi^* c) (\alpha_0, \dots, \alpha_k) := \psi^*(c(\psi^{-1}{}^* \alpha_0, \dots, \psi^{-1}{}^* \alpha_k)).
\]
One can check that this action preserves the differential.

\begin{remark} The complex $C_{\mathrm{def}}(A)$ can be given a structure of differential graded Lie algebra by introducing the classical \emph{Gerstenhaber bracket}. However, we will not need this additional structure. The interested reader may refer to \cite{crainic:def} for further details. \end{remark}

Finally, take a VB-algebroid $(W \Rightarrow E; A \Rightarrow M)$. Deformations of the VB-algebroid structure are controlled by a subcomplex $C_{\mathrm{def,lin}}(W)$ of $C_{\mathrm{def}}(W)$, the \emph{linear deformation complex} of $W$, defined as follows \cite{etv:infinitesimal, lapastina:def}. Let $h_\lambda$ be the homogeneity structure of the vector bundle $W \to A$. We say that a deformation cochain $\tilde c \in C_{\mathrm{def}}(W)$ is \emph{linear} if and only if
\[
h_\lambda^* \tilde c = \tilde c
\]
for every $\lambda > 0$. The linear deformation complex consists, by definition, of linear deformation cochains.

\begin{remark} \label{rmk:lin_def} A linear deformation cochain with $k$ entries is completely determined by its action on $k$ linear sections and $k-1$ linear sections and a core section and the action of its symbol on $k-1$ linear sections and $k-2$ linear sections and a core section. See \cite{lapastina:def} for a proof. \end{remark}

\subsection{VB-groupoids}\label{Sec:deformations2}

In this section we recall from \cite{mackenzie, gracia-saz:vb2, bursztyn:vec} the basic definitions and properties of VB-groupoids that will be useful later. From now on, we will need the notion of a double vector bundle (DVB for short): we refer to \cite{mackenzie} for definitions and basic properties and to our previous paper \cite{lapastina:def} for notations.

\begin{definition}\label{def:VB-gr} A \emph{VB-groupoid} is a \emph{vector bundle in the category of Lie groupoids}, i.e.~a diagram
\begin{equation}\label{eq:VB-gr}
\begin{array}{r}
\xymatrix{
\mathcal V \ar[d] \ar@<0.4ex>[r] \ar@<-0.4ex>[r] & E \ar[d] \\
\mathcal{G} \ar@<0.4ex>[r] \ar@<-0.4ex>[r] & M}
\end{array},
\end{equation}
where $\mathcal V \rightrightarrows E$ and $\mathcal{G} \rightrightarrows M$ are Lie groupoids, $\mathcal V \to \mathcal{G}$ and $E \to M$ are vector bundles and all the vector bundle structure maps (addition, multiplication, projection and zero section) are Lie groupoid maps. We denote $\tilde{\mathsf s}, \tilde{\mathsf t}, \tilde{\mathsf 1}, \tilde{\mathsf m}, \tilde{\mathsf i}$ the structure maps of $\mathcal V$, $\mathsf s, \mathsf t, \mathsf 1, \mathsf m, \mathsf i$ the structure maps of $\mathcal{G}$, $\tilde{\pi}: \mathcal V \to \mathcal{G}$, $\pi: E \to M$ the vector bundle projections and $\tilde 0: \mathcal{G} \to \mathcal V$, $0: M \to E$ the zero sections. The VB-groupoid (\ref{eq:VB-gr}) will be also denoted $(\mathcal V \rightrightarrows E; \mathcal{G} \rightrightarrows M)$. The groupoid $\mathcal V \rightrightarrows E$ is called the \emph{total groupoid}, $\mathcal{G} \rightrightarrows M$ is called the \emph{base groupoid}. We will sometimes say that \emph{$\mathcal V$ is a VB-groupoid over $\mathcal G$}. \end{definition}

\begin{remark} Actually, some of the conditions in the previous definition are redundant and could be omitted. For details and other equivalent definitions of VB-groupoids, see the discussion in \cite{gracia-saz:vb2}. \end{remark} 

\begin{remark}\label{rmk:hom_str} The definition of VB-groupoid can be greatly simplified using the concept of homogeneity structure.
Indeed, consider a diagram of Lie groupoids and vector bundles like (\ref{eq:VB-gr}). From now on, unless otherwise stated, we denote by $h$ the homogeneity structure of $\mathcal V \to \mathcal G$. It can be shown that such a diagram is a VB-groupoid if and only if, for every $\lambda > 0$, $h_\lambda$ is a Lie groupoid automorphism \cite{bursztyn:vec}. \end{remark}

The Lie theory of VB-groupoids is studied in \cite{bursztyn:vec}. Here the authors show that, applying the usual differentation process to the total and the base groupoids of a VB-groupoid, we end up with a VB-algebroid.

From Definition \ref{def:VB-gr} it follows \cite{li-bland} that the map
\[
\pi^R: \mathcal V \to \mathsf s^* E, \quad v \mapsto (\tilde {\mathsf s} (v), \tilde{\pi} (v))
\]
is a surjective submersion. Hence its kernel is a vector bundle $V^R \to \mathcal{G}$, called the \emph{right-vertical subbundle} of $\mathcal V$. Finally, the \emph{right-core} of $\mathcal V$ is $C^R := \mathsf 1^* (V^R)$. Explicitly, $C^R$ is the set of elements of $\mathcal V$ that project on the units of $\mathcal{G}$ and $\mathsf s$-project to the zero section of $M$: this is analogous to the definition of the core of a DVB.

The right-core fits in a short exact sequence of vector bundles over $\mathcal{G}$:
\begin{equation}\label{eq:coreseq}
0 \longrightarrow \mathsf t^* C^R \overset{j^R}{\longrightarrow} \mathcal V \overset{\pi^R}{\longrightarrow} \mathsf s^* E \longrightarrow 0,
\end{equation} 
where $j^R$ is defined by $j^R(c,g) = c \cdot \tilde {\mathsf 0}_g$. A splitting of such a sequence always exists and gives a non-canonical decomposition $\mathcal V \cong \mathsf s^* E \oplus \mathsf t^* C^R$. Additionally, over the submanifold of units of $\mathcal G$ there is a natural splitting, given by $\tilde {\mathsf 1}: E \to \mathcal V$, and, following \cite{gracia-saz:vb2}, we give the following

\begin{definition} A \emph{right-horizontal lift} of the VB-groupoid (\ref{eq:VB-gr}) is a splitting $h: \mathsf s^* E \to \mathcal V$ of (\ref{eq:coreseq}) that satisfies $h(e, \mathsf 1_x) = \tilde {\mathsf 1}_e$ for all $x \in M$, and $e \in E_x$. A \emph{right-decomposition} of (\ref{eq:VB-gr}) is a direct sum decomposition $\mathcal V \cong \mathsf s^\ast E \oplus \mathsf t^\ast C^R$ that comes from a right-horizontal lift. \end{definition}

The existence of right-horizontal lifts can be proved by a partitions of unity argument, hence \emph{every VB-groupoid admits a (non-canonical) right-decomposition} \cite{gracia-saz:vb2}. 

By exchanging the role of the source and the target, one can similarly define a \emph{left-core} $C^L$. The analogue short exact sequence of (\ref{eq:coreseq}) is
\begin{equation}\label{eq:coreseq_2}
0 \longrightarrow \mathsf s^* C^L \longrightarrow \mathcal V \longrightarrow \mathsf t^* E \longrightarrow 0.
\end{equation}
The splittings of (\ref{eq:coreseq_2}) that restrict to the natural splitting over the units are called \emph{left-horizontal lifts}. Moreover, the involution
\[
F: \mathcal V \to \mathcal V, \quad v \mapsto -v^{-1}
\]
induces an isomorphism of vector bundles between the right-core and the left-core.

In the following, we will always consider the right-core: it will be referred to simply as ``core'' and denoted $C$. The \emph{core-anchor} is the vector bundle map $\partial: C \to E$ defined by $\partial c = \tilde {\mathsf t} (c)$.

Now we review some examples of VB-groupoids that will appear later on in the paper.

\begin{examplex} Let $\mathcal{G} \rightrightarrows M$ be a Lie groupoid and $T\mathcal{G} \rightrightarrows TM$ be its tangent (prolongation) groupoid. Then it is easy to check that
\[
\begin{array}{r}
\xymatrix{
T \mathcal{G} \ar[d] \ar@<0.4ex>[r] \ar@<-0.4ex>[r] & TM \ar[d] \\
\mathcal{G} \ar@<0.4ex>[r] \ar@<-0.4ex>[r] & M}
\end{array},
\]
is a VB-groupoid. The core of $(T\mathcal G \rightrightarrows TM; \mathcal G \rightrightarrows M)$ is, by definition, the Lie algebroid $A$ of $\mathcal{G}$ and the core-anchor is the anchor map $\rho: A \to TM$.
\end{examplex}

\begin{examplex} \label{ex:act_grpd} Let $\mathcal G \rightrightarrows M$ be a Lie groupoid and let $\pi: E \to M$ be a representation of $\mathcal G$. Out of these data, we can define the action groupoid $\mathcal G \ltimes E \rightrightarrows E$. As a manifold, $\mathcal G \ltimes E = \mathcal G  \tensor*[_{\mathsf s}]{\times}{_\pi} E$; the structure maps are given by
\begin{equation}\label{action_grpd}
\begin{aligned}
\tilde{\mathsf s} (g, e) & := e, \\ 
\tilde{\mathsf t} (g, e) & := g \cdot e, \\ 
(h, ge) \cdot (g, e) & := (hg, e).
\end{aligned}
\end{equation}
It is easy to prove that
\[
\begin{array}{r}
\xymatrix{
\mathcal{G} \ltimes E \ar[d] \ar@<0.4ex>[r] \ar@<-0.4ex>[r] & E \ar[d] \\
\mathcal{G} \ar@<0.4ex>[r] \ar@<-0.4ex>[r] & M}
\end{array}
\]
is a VB-groupoid. Its core is trivial. Actually every VB-groupoid with trivial core arises in this way, up to a canonical isomorphism.
\end{examplex}

Let us briefly recall how duality works for VB-groupoids. Let $\mathcal V$ be a VB-groupoid like in (\ref{eq:VB-gr}), with core $C$, and let $\mathcal V^* \to \mathcal{G}$ be the dual vector bundle of $\mathcal V \to \mathcal{G}$. Then we can define the \emph{dual VB-groupoid}
\begin{equation}\label{eq:dual}
\begin{array}{r}
\xymatrix{
\mathcal V^* \ar[d] \ar@<0.4ex>[r] \ar@<-0.4ex>[r] & C^* \ar[d] \\
\mathcal{G} \ar@<0.4ex>[r] \ar@<-0.4ex>[r] & M}
\end{array}
\end{equation}
as follows. The source and the target $\check{\mathsf s}, \check {\mathsf t}: \mathcal V^* \rightrightarrows C^*$ are defined by
\[
\begin{aligned}
\langle \check {\mathsf s}(\varphi), c_1 \rangle & := - \langle \varphi, 0_g \cdot c_1^{-1} \rangle, \\
\langle \check {\mathsf t}(\varphi), c_2 \rangle & := \langle \varphi, c_2 \cdot 0_g \rangle
\end{aligned}
\]
for every $g \in \mathcal{G}, \varphi \in \mathcal V^*_g, c_1 \in C_{\mathsf s(g)}, c_2 \in C_{\mathsf t(g)}$, while the multiplication is defined by
\begin{equation}\label{eq:mult}
\langle \varphi_1 \cdot \varphi_2, v_1 \cdot v_2 \rangle := \langle \varphi_1, v_1 \rangle + \langle \varphi_2, v_2 \rangle
\end{equation}
for all $(g_1,g_2) \in \mathcal G^{(2)}$, $(\varphi_1, \varphi_2) \in (\mathcal V^*)^{(2)}$, $\varphi_i \in \mathcal V^*_{g_i}$, $(v_1, v_2) \in \mathcal V^{(2)}$, $v_i \in \mathcal V_{g_i}$. Here and in the following, $\langle \cdot, \cdot \rangle$ denotes the duality pairing. For details and proofs see \cite{mackenzie}.

\begin{examplex} Let $\mathcal G \rightrightarrows M$ be a Lie groupoid, and let $A \Rightarrow M$ be its Lie algebroid. The dual of the tangent VB-groupoid is the \emph{cotangent VB-groupoid}
\[
\begin{array}{r}
\xymatrix{
T^* \mathcal{G} \ar[d] \ar@<0.4ex>[r] \ar@<-0.4ex>[r] & A^* \ar[d] \\
\mathcal{G} \ar@<0.4ex>[r] \ar@<-0.4ex>[r] & M}
\end{array}.
\]
\end{examplex}

Finally, we describe the linear complex and the VB-complex of a VB-groupoid $(\mathcal V \rightrightarrows E, \mathcal{G} \rightrightarrows M)$ \cite{gracia-saz:vb2}. 

We know that the total groupoid comes with its Lie groupoid complex $(C (\mathcal V), \delta)$. It is easy to check that there is an induced vector bundle structure on $\mathcal V^{(k)} \to \mathcal{G}^{(k)}$, so there is a natural subcomplex $C_{\mathrm{lin}}(\mathcal V)$ of $C(\mathcal V)$, whose $k$-cochains are functions on $\mathcal V^{(k)}$ that are linear over $\mathcal{G}^{(k)}$. Inside $C_{\mathrm{lin}}(\mathcal V)$, there is a distinguished subcomplex $C_{\mathrm{proj}}(\mathcal V)$ consisting of \emph{left-projectable linear cochains} \cite{gracia-saz:vb2}. By definition, a linear cochain $f \in C^k_{\mathrm{lin}}(\mathcal V)$ is left-projectable if
\begin{enumerate}
\item $f(\tilde 0_g, v_2, \dots, v_k) = 0$ for every $(\tilde 0_g, v_2, \dots, v_k) \in \mathcal V^{(k)}$;
\item $f(\tilde 0_g \cdot v_1, \dots, v_k) = f(v_1, \dots, v_k)$ for $(v_1, \dots, v_k) \in \mathcal V^{(k)}$ and $g \in \mathcal G$ such that $\tilde {\mathsf t}(v_1) = 0_{\mathsf s(g)}$.
\end{enumerate}

The complexes $C_{\mathrm{lin}}(\mathcal V)$ and $C_{\mathrm{proj}}(\mathcal V)$ are called the \emph{linear complex} and the \emph{VB-complex} of $\mathcal V$, respectively. Their cohomologies are denoted $H_{\mathrm{lin}}(\mathcal V)$ and $H_{\mathrm{proj}}(\mathcal V)$ and called the \emph{linear cohomology} and the \emph{VB-cohomology} of $\mathcal V$.

\begin{remark} We are adopting the terminology of \cite{delhoyo:morita}. In \cite{gracia-saz:vb2} and \cite{crainic:def2}, instead, the VB-complex and the VB-cohomology of $\mathcal V$ are defined to be $C_{\mathrm{proj}}(\mathcal V^*)$ and $H_{\mathrm{proj}}(\mathcal V^*)$, respectively. \end{remark}

Notice that $C (\mathcal G)$ can be identified with the subcomplex of $C (\mathcal V)$ of fiberwise constant cochains. Then, one can check that the product (\ref{eq:prod}) gives $C_{\mathrm{lin}}(\mathcal V)$ and $C_{\mathrm{proj}}(\mathcal V)$ a $C (\mathcal G)$-DG-module structure.

It turns out that the linear and the VB-cohomology are isomorphic:

\begin{lemma} \label{prop:proj;lin} \cite[Lemma 3.1]{cabrera:hom} The inclusion $C_{\mathrm{proj}}(\mathcal V) \hookrightarrow C_{\mathrm{lin}}(\mathcal V)$ induces an isomorphism of $H (\mathcal G)$-modules in cohomology. \end{lemma}

For our purposes, the VB-complex is particularly important because it gives another description of the deformation complex of a Lie groupoid $\mathcal G$. Indeed, we have:

\begin{proposition}\label{prop:cotangent} \cite[Proposition 3.9]{crainic:def2}
There is an isomorphism of $C (\mathcal G)$-DG-modules
\begin{equation}\label{eq:proj}
\phi: C_{\mathrm{def}}(\mathcal G) \rightarrow C_{\mathrm{proj}} (T^* \mathcal G)[1]
\end{equation}
given by
\[
\phi(c)(\theta_0, \dots, \theta_k) = \langle \theta_0, c(g_0, \dots, g_k) \rangle
\]
for all $c \in C^k_{\mathrm{def}}(\mathcal G)$ and $(\theta_0, \dots, \theta_k) \in (T^* \mathcal G)^{(k+1)}$ such that $\theta_i \in T^*_{g_i} \mathcal G$.
\end{proposition}

\subsection{The linear deformation complex of a VB-groupoid}\label{sec:lin_def}

In this subsection we introduce the main object of this paper: the \emph{linear deformation complex of a VB-groupoid}, first introduced in \cite{etv:infinitesimal} (for different purposes from the present ones). The definition is entirely analogous to the one recalled in Subsection \ref{Sec:deformations1} for VB-algebroids.
 
Let $(\mathcal V \rightrightarrows E; \mathcal G \rightrightarrows M)$ be a VB-groupoid, and let $(W \Rightarrow E; A \Rightarrow M)$ be its VB-algebroid. By Remark \ref{rmk:hom_str}, $h_\lambda$ is a Lie groupoid automorphism for every $\lambda > 0$, so it acts on the deformation complex $C_{\mathrm{def}}(\mathcal V)$ of $\mathcal V \rightrightarrows E$. We say that a deformation cochain $\tilde c$ is \emph{linear} if
\begin{equation}\label{eq:lin_coch}
h_\lambda^* \tilde c = \tilde c
\end{equation}
for every $\lambda > 0$. Hence, linear cochains are those which are invariant under the homogeneity structure.

We know that $h_\lambda^*$ commutes with $\delta$, for all $\lambda > 0$, so \emph{linear deformation cochains form a subcomplex of $C_{\mathrm{def}}(\mathcal V)$}. We denote the latter by $C_{\mathrm{def,lin}}(\mathcal V)$ and we call it the \emph{linear deformation complex} of $\mathcal V$. Its cohomology is called the \emph{linear deformation cohomology} of $\mathcal V$ and denoted $H_{\mathrm{def,lin}}(\mathcal V)$. Formula (\ref{eq:prod}) also shows that \emph{$C_{\mathrm{def,lin}}(\mathcal V)$ is a $C (\mathcal G)$-module}.

The action of $h_\lambda$ on $C^{-1}_{\mathrm{def}}(\mathcal V)$ coincides, by definition, with the action induced by the homogeneity structure of $W \to A$ on $\Gamma(W,E)$, so $C^{-1}_{\mathrm{def,lin}}(\mathcal V)$ is simply the space $\Gamma_{\mathrm{lin}}(W,E)$ of linear sections of $W \to E$. For $k \geq 0$, Equation (\ref{eq:lin_coch}) is equivalent to saying that $\tilde c: \mathcal V^{(k+1)} \to T \mathcal V$ intertwines the homogeneity structures of $\mathcal V^{(k+1)} \to \mathcal G^{(k+1)}$ and $T \mathcal V \to T \mathcal G$. Again by Remark \ref{rmk:hom_str}, this means that $\tilde c$ is a vector bundle map over some map $c: \mathcal G^{(k+1)} \to T \mathcal G$. In this way we recover the definition in \cite{etv:infinitesimal}. 

For $k \geq 0$, a linear $k$-cochain $\tilde c$ can also be seen as an $\tilde {\mathsf s}$-projectable, linear section of the DVB $(p^*_{k+1} T \mathcal V \to \mathcal V^{(k+1)}; p^*_{k+1} T \mathcal G  \to \mathcal G^{(k+1)})$:
\begin{equation}\label{eq:pullback_DVB}
\begin{array}{c}
\xymatrix{\tilde p^*_{k+1} T \mathcal V \ar[r] \ar[d] & \mathcal V^{(k+1)} \ar[d] \ar@/_1.25pc/[l]_-{\tilde c}\\
p^*_{k+1} T \mathcal G \ar[r] & \mathcal G^{(k+1)} \ar@/_1.25pc/[l]_-{ c}}
\end{array}
,
\end{equation}
where we denote $\tilde p_k: \mathcal V^{(k)} \to \mathcal V$, and $\tilde q_k: \mathcal V^{(k)} \to E$ the maps (\ref{eq:pq}) for $\mathcal V$.

There is another way to describe the linear deformation complex of a VB-groupoid. Let $(\mathcal V \rightrightarrows E; \mathcal G \rightrightarrows M)$ be a VB-groupoid, let $C$ be its core, and let $(W \Rightarrow E; A \Rightarrow M)$ be its VB-algebroid. Consider the cotangent VB-groupoid of $\mathcal V \rightrightarrows E$, $(T^* \mathcal V \rightrightarrows W^*_E; \mathcal V \rightrightarrows E)$ (here we denote by $W^\ast_E \to E$ the dual of a vector bundle $W \to E$). 

Actually, $T^* \mathcal V \rightrightarrows W^*_E$ is the top groupoid of another VB-groupoid. To see this, first take the dual of $\mathcal V$, that is $(\mathcal V^* \rightrightarrows C^*; \mathcal G \rightrightarrows M)$. Then, the cotangent VB-groupoid of $\mathcal V^* \rightrightarrows C^*$ is $(T^* \mathcal V^* \rightrightarrows (W^*_A)^*_{C^*}; \mathcal V^* \rightrightarrows C^*)$. Finally, recall from \cite{mackenzie} that there is a canonical isomorphism 
\[ 
B: T^*\mathcal V \to T^* \mathcal V^*
\]
of both DVBs and Lie groupoids, covering an isomorphism
\[
\beta : W^\ast_E \to (W^\ast_{C^\ast})^\ast_{C^\ast}.
\]
of DVBs. Combining all these maps, we obtain a diagram:

\begin{equation}\label{eq:dvb-grpd}
\begin{array}{c}
\xymatrix@C=10pt@R=15pt{
 & T^\ast \mathcal V \ar[dd]|!{[dl];[dr]}{\hole} \ar@<0.4ex>[rr] \ar@<-0.4ex>[rr] \ar[dl]^-{B}_-{\cong}& & W^\ast_E \ar[dl]^-{\beta}_-{\cong} \ar[dd] \\
 T^\ast \mathcal V^\ast \ar[dd] \ar@<0.4ex>[rr] \ar@<-0.4ex>[rr] & & (W^*_A)^*_{C^*} \ar[dd]& \\
 & V^\ast \ar@<0.4ex>[rr]|!{[ur];[dr]}{\hole} \ar@<-0.4ex>[rr]|!{[ur];[dr]}{\hole} \ar@{=}[dl]& & C^\ast \ar@{=}[dl] \\
 V^\ast \ar@<0.4ex>[rr] \ar@<-0.4ex>[rr] & & C^* & 
}
\end{array}.
\end{equation}

The maps $B$ and $\beta$ are isomorphisms of vector bundles over the identity, so the back face in the diagram (\ref{eq:dvb-grpd}) is also a VB-groupoid.

It follows that inside $C (T^* \mathcal V)$ there are two distinguished subcomplexes, those of cochains that are linear over $\mathcal V$ and over $\mathcal V^*$: we denote them by $C_{\mathrm{lin},\bullet}(T^* \mathcal V)$ and $C_{\bullet, \mathrm{lin}}(T^* \mathcal V)$ respectively. Moreover, denote $C_{\mathrm{proj},\bullet}(T^* \mathcal V)$ the subcomplex of left-projectable linear cochains over $\mathcal V$ and define $C_{\mathrm{lin,lin}}(T^* \mathcal V) := C_{\mathrm{lin}, \bullet}(T^* \mathcal V) \cap C_{\bullet, \mathrm{lin}}(T^* \mathcal V)$, $C_{\mathrm{proj,lin}}(T^* \mathcal V) := C_{\mathrm{proj}, \bullet}(T^* \mathcal V) \cap C_{\bullet, \mathrm{lin}}(T^* \mathcal V)$. We denote their cohomologies $H_{\mathrm{lin,lin}}(\mathcal V)$ and $H_{\mathrm{proj,lin}}(\mathcal V)$ respectively.

From Proposition \ref{prop:cotangent}, there is an isomorphism of $C (\mathcal V)$-modules
\begin{equation}\label{eq:iso}
C_{\mathrm{def}}(\mathcal V) \cong C_{\mathrm{proj},\bullet} (T^* \mathcal V)[1].
\end{equation}

It is easy to check that this isomorphism takes linear deformation cochains to cochains on $T^* \mathcal V$ that are linear over $\mathcal V^*$. So we get the following

\begin{proposition}\label{prop:cotangent_VB} There is an isomorphism of $C (\mathcal G)$-modules
\begin{equation}\label{eq:iso_lin}
C_{\mathrm{def,lin}}(\mathcal V) \cong C_{\mathrm{proj,lin}}(T^* \mathcal V)[1].
\end{equation}
\end{proposition}

For later use, we notice that a ``linear version'' of Lemma \ref{prop:proj;lin} holds. Namely, we have

\begin{lemma}\label{prop:proj,lin;lin,lin} The inclusion $C_{\mathrm{proj,lin}}(T^* \mathcal V) \hookrightarrow C_{\mathrm{lin,lin}}(T^* \mathcal V)$ induces an isomorphism in cohomology. \end{lemma}

\proof The proof of \cite[Lemma 3.1]{cabrera:hom} works identically in our setting without significant modifications. \endproof

\subsubsection{Deformations of $\mathcal G$ from linear deformations of $\mathcal V$} 
We have the following

\begin{proposition}\cite[Lemma 2.27]{etv:infinitesimal} If $\tilde c: \mathcal V^{(k+1)} \to T\mathcal V$ belongs to $C^k_{\mathrm{def,lin}}(\mathcal V)$, then its projection $c: \mathcal{G}^{(k+1)} \to T\mathcal{G}$ belongs to $C^k_{\mathrm{def}}(\mathcal{G})$ and $\delta \tilde c$ projects to $\delta c$. \end{proposition}

It follows that there exists a natural cochain map:
\begin{equation}\label{eq:VB_pr}
C_{\mathrm{def,lin}}(\mathcal V) \to C_{\mathrm{def}}(\mathcal{G}).
\end{equation}
In degree $k = -1$, this is simply the projection $\Gamma_{\mathrm{lin}}(W,E) \to \Gamma(A)$ and we have the well-known short exact sequence:
\begin{equation} \label{eq:lin_sec}
0 \longrightarrow \mathfrak{Hom} (E,C) \longrightarrow \Gamma_{\mathrm{lin}}(W,E) \longrightarrow \Gamma(A) \longrightarrow 0,
\end{equation}
where $\mathfrak{Hom} (E,C)$ is the $C^\infty(M)$-module of vector bundle morphisms $E \to C$. 

We now show that the map (\ref{eq:VB_pr}) is surjective for all $k \geq 0$. Let $c \in C^k_{\mathrm{def}}(\mathcal G)$. By Remark \ref{rmk:couple}, we have a diagram:
\[
\begin{array}{c}
\xymatrix{\mathcal G^{(k+1)} \ar[r]^{c} \ar[d] & p_{k+1}^* T \mathcal G \ar[d]^{T\mathsf s} \\
\mathcal G^{(k)} \ar[r]^{s_c} & q_k^* TM}
\end{array}
.
\]
We can lift $s_c$ to a linear section $s_{\tilde c}$ of $\tilde q_k^* TE \to \mathcal V^{(k)}$, so we obtain the following diagram
\[
\begin{array}{c}
\xymatrix@C=9pt@R=12pt{
 & \tilde p^*_{k+1} T \mathcal V \ar[dd]|!{[dl];[dr]}{\hole} \ar[rr] \ar[dl] & & \tilde q_k^* TE \ar[dl] \ar[dd] \\
 \mathcal V^{(k+1)} \ar[dd] \ar[rr]  & & \mathcal V^{(k)} \ar@/_0.8pc/[ur]_-{s_{\tilde c}} \ar[dd]& \\
 & p^*_{k+1} T \mathcal G \ar[rr]|!{[ur];[dr]}{\hole} \ar[dl]& & q_k^* TM \ar[dl] \\
 \mathcal G^{(k+1)} \ar[rr]  \ar@/_0.8pc/[ur]_-{c}& & \mathcal G^{(k)} \ar@/_0.8pc/[ur]_-{s_{ c}}& 
}
\end{array},
\]
where the left and the right faces are DVBs and the horizontal arrows form a surjective DVB morphism. Thus we are in the situation of Lemma \ref{lemma:section} and we conclude that there exists a linear section $\tilde c: \mathcal V^{(k+1)} \to \tilde p^*_{k+1} T \mathcal V$ that projects on $c$, as desired. Summarizing, there is a canonical short exact sequence of cochain maps:
\begin{equation}\label{eq:ses_VB}
0 \longrightarrow \ker \Pi \longrightarrow C_{\mathrm{def,lin}}(\mathcal V) \overset{\Pi}{\longrightarrow} C_{\mathrm{def}} (\mathcal{G})\longrightarrow 0.
\end{equation}

Now we compute $\ker \Pi$. By definition, a $k$-cochain $\tilde c: \mathcal V^{(k+1)} \to T\mathcal V$ is killed by $\Pi$ if and only if it takes values in the vertical bundle $T^{\tilde \pi} \mathcal V$ of $\tilde \pi: \mathcal V \to \mathcal G$. It is easy to check that $T^{\tilde \pi} \mathcal V \rightrightarrows T^\pi E$ is a subgroupoid of $T\mathcal V \rightrightarrows TE$. Moreover, $T^{\tilde \pi} \mathcal V \cong \mathcal V \times_{\mathcal G} \mathcal V$ and $T^\pi E \cong E \times_M E$ canonically as vector bundles. Under these isomorphisms, $T^{\tilde \pi} \mathcal V$ is identified with the groupoid $\mathcal V \times_{\mathcal G} \mathcal V \rightrightarrows E \times_M E$ with the component-wise structure maps. We will understand this identification.

It is clear that $\tilde c (v_0, \dots, v_k)$ has $v_0$ as first component, so one can think of $\tilde c$ as a map $\mathcal V^{(k+1)} \to \mathcal V$ that is linear over $p_{k+1}: \mathcal G^{(k+1)} \to \mathcal G$. Then $\ker \Pi$ is given by elements $\tilde c \in \mathfrak{Hom} (\mathcal V^{(k+1)}, p_{k+1}^* \mathcal V)$ such that $\tilde {\mathsf s} (c(v_0, \dots, v_k))$ does not depend on $v_0$ for any $(v_0, \dots, v_k) \in \mathcal V^{(k+1)}$. Finally, we observe that, on $\ker \Pi$, the differential is simply given by
\[
\begin{aligned}
\delta \tilde c (v_0, \dots, v_{k+1})  = &-\tilde{\bar{\mathsf m}} (\tilde c(v_0 v_1, \dots, v_{k+1}), \tilde c(v_1, \dots, v_{k+1})) \\ 
& + \sum_{i=1}^k (-1)^{i-1} \tilde c(v_0, \dots, v_i v_{i+1}, \dots, v_{k+1}) + (-1)^{k} \tilde c(v_0, \dots, v_k).
\end{aligned}
\]

\subsubsection{Trivial-core VB-groupoids}

A VB-groupoid $(\mathcal V \rightrightarrows E; \mathcal G \rightrightarrows M)$ is called \emph{trivial-core} (or \emph{vacant}) if its core is the zero-vector bundle $0_M \to M$. We want to show that the linear deformation complex of a trivial-core VB-groupoid has a particularly simple shape. First recall that, in Example \ref{ex:act_grpd}, we have observed that the total groupoid of a trivial-core VB-groupoid is canonically isomorphic to the action groupoid $\mathcal G \ltimes E \rightrightarrows E$  associated to a representation of the base groupoid $\mathcal G \rightrightarrows M$ on the side bundle $E$. Therefore, without loss of generality, we will assume that $\mathcal V = \mathcal G \ltimes E$. As a vector bundle over $\mathcal G$, it is the pull-back $\mathsf s^\ast E$, and we denote it also by $E_{\mathcal G}$.

Consider a linear cochain $\tilde c \in C^k_{\mathrm{def,lin}}(E_{\mathcal G})$. By definition, it gives a commutative diagram:
\begin{equation}\label{diag:c_tilde}
\begin{array}{c}
\xymatrix@C=7pt@R=7pt{
& E_{\mathcal G}^{(k+1)} \ar[dl] \ar[rr]^{\tilde c} \ar[dd]|!{[dl];[dr]}{\hole} & & TE_{\mathcal G} \ar[dl] \ar[dd] \\
E_{\mathcal G}^{(k)} \ar[rr]^(.6){s_{\tilde c}} \ar[dd] & & TE \ar[dd] \\
& \mathcal G^{(k+1)} \ar[dl] \ar[rr]^(.4)c |!{[ur];[dr]}{\hole} & & T \mathcal G  \ar[dl] \\
\mathcal G^{(k)} \ar[rr]^(.6){s_c} & & TM  }
\end{array}.
\end{equation}
But, for every $k$, there is a canonical isomorphism 
\begin{equation} \label{eq:nerve_act}
E_{\mathcal G}^{(k)} \overset{\cong}{\longrightarrow} \mathcal G^{(k)} \times_M E, \quad ((g_1, e_1), \ldots, (g_k, e_k)) \longmapsto ((g_1, \ldots, g_k); e_k)
\end{equation}
where the fibered product is wrt the projection $\mathcal G^{(k)} \to M$, $(g_1, \ldots, g_k) \mapsto \mathsf s(g_k)$. We also have that $ TE_{\mathcal G} \cong T \mathcal G \times_{TM} TE$ where the fibered product is wrt to $T \mathsf s : TG \to TM$. So, we get the following alternative description of (\ref{diag:c_tilde}):
\begin{equation} \label{diag:c_tilde_2}
\begin{array}{c}
\xymatrix@C=0pt@R=7pt{
& \mathcal G^{(k+1)} \times_M E \ar[dl] \ar[rr]^{\tilde c} \ar[dd]|!{[dl];[dr]}{\hole} & & T \mathcal G \times_{TM} TE \ar[dl] \ar[dd] \\
\mathcal G^{(k)} \times_M E \ar[rr]^(.65){s_{\tilde c}} \ar[dd] & & TE \ar[dd] \\
& \mathcal G^{(k+1)} \ar[dl] \ar[rr]^(.4)c |!{[ur];[dr]}{\hole} & & T \mathcal G  \ar[dl] \\
\mathcal G^{(k)} \ar[rr]^(.65){s_c} & & TM  }
\end{array},
\end{equation}
where the vertical arrows, except for the front right one, are projections onto the first factor. In particular, $\tilde c$ is fully determined by $c$ and $s_{\tilde c}$. Set $\tilde c_1 := c$ and observe that, for every $(g_1, \dots, g_k) \in \mathcal G^{(k)}$, $s_{\tilde c} ((g_1, \dots, g_k); -)$ is a linear map $E_{\mathsf s (g_k)} \to TE|_{s_c (g_1, \dots, g_k)}$. Consider the linear map
\begin{equation} \label{eq:c_2}
\tilde c_2 (g_1, \dots, g_k): E_{\mathsf t (g_1)} \to TE|_{s_c (g_1, \dots, g_k)}, \quad e \mapsto s_{\tilde c} ((g_1, \dots, g_k); (g_1 \cdots g_k)^{-1} \cdot e).
\end{equation}
It is easy to see that $\tilde c_2 (g_1, \dots, g_k)$ splits the projection $ TE|_{s_c (g_1, \dots, g_k)} \to E_{\mathsf t (g_1)}$. Hence, the pair $(s_c (g_1, \dots, g_k), \tilde c_2 (g_1, \dots, g_k))$ corresponds to a derivation in $D_{\mathsf t (g_1)} E$ (via the inverse to the correspondence $\delta \mapsto (\sigma_\delta, \widehat \delta)$). We denote by $\tilde c_2 (g_1, \dots, g_k)$ again the latter derivation. In this way, we have defined a map $\tilde c_2: \mathcal G^{(k)} \to DE$ such that
\begin{enumerate}
		\item[(TC1)] $\tilde c_2 (g_1, \dots, g_k) \in D_{\mathsf t (g_1)} E$ for every $(g_1, \dots, g_k) \in \mathcal G^{(k)}$; 
		\item[(TC2)] $\sigma \circ \tilde c_2 = s_{\tilde c_1}$.
\end{enumerate}
Conversely, from a pair $(\tilde c_1, \tilde c_2)$ with $\tilde c_1 \in C^k_{\mathrm{def}}(\mathcal G)$ and $\tilde c_2 : \mathcal G^{(k)} \to DE$ satisfying (TC1) and (TC2) above, we can reconstruct a linear deformation cochain $\tilde c \in C^k_{\mathrm{def}, \mathrm{lin}} (E_{\mathcal G})$ in the obvious way. Finally, a direct computation exploiting \cite[Formula (2)]{crainic:def2} shows that, for every $\tilde c \in C^k_{\mathrm{def}, \mathrm{lin}} (E_{\mathcal G})$,
\begin{equation}\label{eq:delta_c_1}
(\delta \tilde c)_1 = \delta (\tilde c_1)
\end{equation}
and $(\delta \tilde c)_2$ is given by the following formula
\begin{equation}\label{eq:delta_c_2}
\begin{aligned}
& (\delta \tilde c)_2 (g_1, \dots, g_{k+1}) (e) \\
& =  - \tilde c_1(g_1, \dots, g_{k+1}) * \left(\tilde c_2 (g_2, \dots, g_{k+1}) (g_1^{-1} e)\right) \\
& \quad + \sum_{i=1}^k (-1)^{i-1} \tilde c_2 (g_1, \dots, g_i g_{i+1}, \dots, g_{k+1})(e) + (-1)^k \tilde c_2 (g_1, \dots, g_k)(e),
\end{aligned}
\end{equation}
where the ``$\ast$'' is the tangent map $\ast : T \mathcal G \times_{TM} TE \to TE$ to the action $\tilde{\mathsf t} : \mathcal G \ltimes E = \mathcal G \times_M E \to E$.

The above discussion proves the following

\begin{lemma} \label{lem:trivial_core} Let $(\mathcal G \ltimes E \rightrightarrows E; \mathcal G \rightrightarrows E)$ be a trivial-core VB-groupoid and let $k \geq 0$. The assignment $\tilde c \mapsto (\tilde c_1, \tilde c_2)$ establishes an isomorphism between $C^k_{\mathrm{def,lin}}(E_{\mathcal G})$ and the space of pairs $(\tilde c_1, \tilde c_2)$ with $\tilde c_1 \in C_{\mathrm{def}}^k(\mathcal G)$, and $\tilde c_2 : \mathcal G^{(k)} \to DE$ satisfying \emph{(TC1)} and \emph{(TC2)} above. Under this isomorphism the differential $\delta \tilde c$ corresponds to the pair $((\delta \tilde c)_1, (\delta \tilde c)_2)$ given by Formulas (\ref{eq:delta_c_1}) and (\ref{eq:delta_c_2}).
\end{lemma}

Now, take $\tilde c \in C^k_{\mathrm{def,lin}}(E_{\mathcal G})$. Using that $DE_{\mathcal G} \cong T \mathcal G \times_{TM} DE$ define a map $\widehat{\tilde c} : \mathcal G^{(k+1)} \to DE_{\mathcal G}$ by putting
\[
\widehat{\tilde c} (g_0, \ldots, g_k) := \left(\tilde c_1 (g_0, \ldots, g_k), \tilde c_2 (g_1, \ldots, g_k) \right),
\]
and observe that 
\begin{enumerate}
\item[(TC3)] $\widehat{\tilde c} (g_0, \dots, g_k) \in D_{g_0} E_{\mathcal G}$, for all $(g_0, \ldots, g_k) \in \mathcal G^{(k+1)}$;
\item[(TC4)] there exists a (necessarily unique) smooth map $\mathcal G^{(k)} \to DE$ making the following diagram commutative:
	\[
\begin{array}{c}
\xymatrix{\mathcal G^{(k+1)} \ar[r]^{\widehat{\tilde c}} \ar[d] & DE_{\mathcal G} \ar[d]^{D \tilde {\mathsf s}} \\
\mathcal G^{(k)} \ar[r] & DE}
\end{array}
,
\]
where the map on the left is the projection onto the last $k$ arrows.
\end{enumerate}
Conversely, given a map $\widehat{\tilde c} : \mathcal G^{(k+1)} \to DE_{\mathcal G}$ satisfying (TC3) and (TC4) we can reconstruct $\tilde c_1, \tilde c_2$ and hence $\tilde c$. A direct computation exploiting Formulas (\ref{eq:delta_c_1}) and (\ref{eq:delta_c_1}) shows that $\widehat{\delta \tilde c} : \mathcal G^{(k+2)} \to DE_{\mathcal G}$ is given by
\begin{equation}\label{eq:delta_c_hat}
\begin{aligned}
\widehat{\delta \tilde c}  (g_0, \dots, g_{k+1}) = & -D \tilde{\bar{\mathsf m}} \left(\widehat{\tilde c}(g_0 g_1, \dots, g_{k+1}), \widehat{\tilde c}(g_1, \dots, g_{k+1})\right) + \\ 
& \sum_{i=1}^k (-1)^{i-1} \widehat{\tilde c}(g_0, \dots, g_i g_{i+1}, \dots, g_{k+1}) + (-1)^k \widehat{\tilde c}(g_0, \dots, g_k),
\end{aligned}
\end{equation}
where we used that $\tilde{\bar{\mathsf m}} : E_{\mathcal G} \mathbin{{}_{\tilde{\mathsf s}} \times_{\tilde{\mathsf s}}} E_{\mathcal G} \to E_{\mathcal G}$ is a regular vector bundle map (covering $\bar{\mathsf m} : \mathcal G \mathbin{{}_{\mathsf s} \times_{\mathsf s}} \mathcal G \to \mathcal G$) and $D(E_{\mathcal G} \mathbin{{}_{\tilde{\mathsf s}} \times_{\tilde{\mathsf s}}} E_{\mathcal G}) \cong DE_{\mathcal G} \mathbin{{}_{D\tilde{\mathsf s}} \times_{D\tilde{\mathsf s}}}DE_{\mathcal G}$ (we leave the easy details to the reader). Summarizing, we have the following

\begin{corollary} \label{cor:trivial_core}
Let $(\mathcal G \ltimes E \rightrightarrows E; \mathcal G \rightrightarrows E)$ be a trivial-core VB-groupoid and let $k \geq 0$. The assignment $\tilde c \mapsto \widehat{\tilde c}$ establishes an isomorphism between $C^k_{\mathrm{def,lin}}(E_{\mathcal G})$ and the space of maps $\widehat{\tilde c} : \mathcal G^{(k+1)} \to DE_{\mathcal G}$ satisfying \emph{(TC3)} and \emph{(TC4)}. Under this isomorphism the differential $\delta \tilde c$ corresponds to the map $\widehat{\delta \tilde c} : \mathcal G^{(k+2)} \to DE_{\mathcal G}$ given by Formula (\ref{eq:delta_c_hat}).
\end{corollary}
Notice the analogy between (\ref{eq:delta_c_hat}) and (\ref{eq:diff}).

Finally, in this case, the projection (\ref{eq:VB_pr}) is the map $\tilde c \mapsto \tilde c_1$. From the condition (TC2), a cochain in the kernel of this projection is equivalent to a map $\tilde c_2: \mathcal G^{(k)} \to DE$ such that $\sigma \circ \tilde c_2 = 0$, so $\tilde c_2$ takes values in $\operatorname{End} E$, and it is easy to see that the sequence (\ref{eq:ses_VB}) takes the form:
\begin{equation} \label{eq:coch_2}
0 \longrightarrow C(\mathcal G, \operatorname{End} E) \longrightarrow C_{\mathrm{def,lin}}(\mathcal G \ltimes E) \longrightarrow C_{\mathrm{def}}(\mathcal G) \to 0,
\end{equation}
where $C(\mathcal G, \operatorname{End} E)$ is the Lie groupoid complex of $\mathcal G$ with ceofficients in the representation of $\mathcal G$ on $\operatorname{End} E$ induced by that on $E$.

\subsubsection{Low-degree cohomology groups}

Next we describe low-degree cohomology groups. The entire discussion of Subsection \ref{Sec:deformations1} goes through, with minor changes. We report it here for completeness.

Let $(\mathcal V \rightrightarrows E; \mathcal G \rightrightarrows M)$ be a VB-groupoid, and let $(W \Rightarrow E; A \Rightarrow M)$ be its VB-algebroid. Consider the isotropy bundle $\mathfrak i$ of $\mathcal V$ and its sections. It is natural to define \emph{linear sections} of $\mathfrak i$:
\[
\Gamma_{\mathrm{lin}}(\mathfrak i) := \Gamma(\mathfrak i) \cap \Gamma_{\mathrm{lin}}(W,E)
\]
and \emph{linear invariant sections}:
\[
H^0_{\mathrm{lin}}(\mathcal V, \mathfrak i) = \Gamma_{\mathrm{lin}}(\mathfrak i)^{\mathrm{inv}} := \Gamma(\mathfrak i)^{\mathrm{inv}} \cap \Gamma_{\mathrm{lin}}(W,E).
\]
The following proposition is easily proved as in \cite[Proposition 4.1]{crainic:def2}.

\begin{proposition}\label{prop:H^-1_lin} Let $(\mathcal V \rightrightarrows E; \mathcal G \rightrightarrows M)$ be a VB-groupoid. Then $H^{-1}_{\mathrm{def,lin}}(\mathcal V) \cong H^0_{\mathrm{lin}}(\mathcal V, \mathfrak i) = \Gamma_{\mathrm{lin}}(\mathfrak i)^{\mathrm{inv}}$. \end{proposition}

Now we consider the complex $C (\mathcal V, \mathfrak i)$. It is clear that there is a distinguished subcomplex $C_{\mathrm{lin}}(\mathcal V, \mathfrak i)$, the one of cochains $\mathcal V^{(k)} \to W$ that are linear over some map $\mathcal G^{(k)} \to A$. From a direct computation it follows that the map (\ref{eq:isotropy}) takes $C_{\mathrm{lin}}(\mathcal V, \mathfrak i)$ to $C_{\mathrm{def,lin}}(\mathcal V)$, so we have a map
\[
r: C_{\mathrm{lin}}(\mathcal V, \mathfrak i) \hookrightarrow C_{\mathrm{def,lin}}(\mathcal V).
\]

\begin{proposition}[see also {\cite[Section 2.3]{etv:infinitesimal}}]
\[
H^0_{\mathrm{def,lin}}(\mathcal V) = \dfrac{\text{linear multiplicative vector fields on $\mathcal V$}}{\text{inner linear multiplicative vector fields on $\mathcal V$}}.
\]
\end{proposition}

\proof Let $X \in C^0_{\mathrm{def,lin}}(\mathcal V)$. Then $X$ is a linear vector field, and from Proposition \ref{prop:H^0} it follows that $X$ is closed if and only if it is multiplicative, and is exact if and only if is inner multiplicative, as desired. \endproof

Recall that also the normal bundle $\nu$ of $\mathcal V$ is defined. Observe that the anchor $\rho: W \to TE$ of the Lie algebroid $W$ is a morphism of DVBs, hence it takes linear sections to linear vector fields $\mathfrak X_{\mathrm{lin}} (E)$. We set
\[
\Gamma_{\mathrm{lin}}(\nu) := \mathfrak X_{\mathrm{lin}}(E)/ \rho(\Gamma_{\mathrm{lin}}(W,E)).
\]

Following the discussion in Subsection \ref{Sec:deformations1}, we declare that a section $[V] \in \Gamma_{\mathrm{lin}}(V)$ is \emph{invariant} if it possesses an $(\mathsf s, \mathsf t)$-lift $X \in \mathfrak X_{\mathrm{lin}}(\mathcal V)$. The space of invariant linear sections is denoted $H^0_{\mathrm{lin}} (\mathcal V, \nu)$ or $\Gamma_{\mathrm{lin}}(\nu)^{\mathrm{inv}}$. Observing that the projection on $E$ of a linear multiplicative vector field is linear, we obtain a linear map
\[
\pi: H^0_{\mathrm{def,lin}}(\mathcal V) \to \Gamma_{\mathrm{lin}}(\nu)^{\mathrm{inv}}.
\]

From Lemma \ref{prop:curvature}, Proposition \ref{prop:exact} and their proofs in \cite{crainic:def2}, the ``linear versions'' follow immediately.

\begin{lemma} Let $[V] \in \Gamma_{\mathrm{lin}}(\nu)^{\mathrm{inv}}$ and $X \in \mathfrak X_{\mathrm{lin}}(\mathcal V)$ an $(\mathsf s, \mathsf t)$-lift of $\mathcal V$. Then $\delta X \in C^2_{\mathrm{lin}}(\mathcal V, \mathfrak i)$ and its cohomology class does not depend on the choice of $X$, hence there is an induced linear map
\[
K: \Gamma_{\mathrm{lin}}(\nu)^{\mathrm{inv}} \to H^2_{\mathrm{lin}}(\mathcal V, \mathfrak i).
\]
\end{lemma}

\begin{proposition} There is an exact sequence
\begin{equation}\label{eq:exact_lin} 
0 \longrightarrow H^1_{\mathrm{lin}}(\mathcal V, \mathfrak i) \overset{r}{\longrightarrow} H^0_{\mathrm{def,lin}}(\mathcal V) \overset{\pi}{\longrightarrow} \Gamma_{\mathrm{lin}}(\nu)^{\mathrm{inv}} \overset{K}{\longrightarrow} H^2_{\mathrm{lin}}(\mathcal V, \mathfrak i) \overset{r}{\longrightarrow} H^1_{\mathrm{def,lin}}(\mathcal V).
\end{equation}
\end{proposition}

\subsubsection{Deformations}

Here we describe deformations of a VB-groupoid $(\mathcal V \rightrightarrows E; \mathcal G \rightrightarrows M)$ and their relation with the $1$-cohomology $H^1_{\mathrm{def,lin}}(\mathcal V)$. Let $B$ be a smooth manifold.

\begin{definition} A \emph{family of VB-groupoids} over $B$ is a diagram
\[
\begin{array}{r}
\xymatrix{
\tilde{\mathcal V} \ar[d] \ar@<0.4ex>[r] \ar@<-0.4ex>[r] & \tilde E \ar[d] \ar[r] & B \ar@{=}[d] \\
\tilde{\mathcal{G}} \ar@<0.4ex>[r] \ar@<-0.4ex>[r] & \tilde M \ar[r] & B }
\end{array}
\]
such that the first square is a VB-groupoid and the rows are families of Lie groupoids. In particular, for every $b \in B$, 
\[
\begin{array}{r}
\xymatrix{
\mathcal V_b \ar[d] \ar@<0.4ex>[r] \ar@<-0.4ex>[r] & E_{b} \ar[d] \\
\mathcal{G}_{b} \ar@<0.4ex>[r] \ar@<-0.4ex>[r] & M_{b}}
\end{array}
\]
is a VB-groupoid.

If $B$ is an open interval $I$ containing 0, the family is said to be a \emph{deformation} of $\mathcal V_0$ and the latter is denoted simply by $(\mathcal V \rightrightarrows E; \mathcal G \rightrightarrows M)$. A deformation of $\mathcal V$ is also denoted $(\mathcal V_{\epsilon})$.
\end{definition}

The structure maps of $\mathcal V_{\epsilon}$ are denoted $\tilde {\mathsf s}_{\epsilon}, \tilde {\mathsf t}_{\epsilon}, \tilde {\mathsf 1}_{\epsilon}, \tilde {\mathsf m}_{\epsilon}, \tilde {\mathsf i}_{\epsilon}$, the division map is denoted $\tilde{\bar {\mathsf m}}_{\epsilon}$. Strict, $\mathsf s$-constant, $\mathsf t$-constant, $(\mathsf s ,\mathsf t)$-constant and constant deformations are defined as in the plain groupoid case.

Two deformations $(\mathcal V_{\epsilon})$ and $(\mathcal V'_{\epsilon})$ of $\mathcal V$ are called \emph{equivalent} if there exists a smooth family of VB-groupoid isomorphisms $\Psi_{\epsilon}: \mathcal V_{\epsilon} \to \mathcal V'_{\epsilon}$ such that $\Psi_0 = \mathrm{id}$. We say that $(\mathcal V_{\epsilon})$ is \emph{trivial} if it is equivalent to the constant deformation.

\begin{proposition} [$(\mathsf s, \mathsf t)$-constant deformations] Let $(\mathcal V_{\epsilon})$ be an $(\mathsf s, \mathsf t)$-constant deformation of the VB-groupoid $\mathcal V$. Then formula
\[
u_0 (v, v') := - \dfrac{d}{d \epsilon} \bigg|_{\epsilon = 0} R^{-1}_{v v'}(\tilde {\mathsf m}_{\epsilon}(v,v'))
\]
defines a cocycle $u_0 \in C^1_{\mathrm{lin}}(\mathcal V, \mathfrak i)$ (here $R$ denotes right translation in $\mathcal V$). Its image $\xi_0$ in $C^1_{\mathrm{def,lin}}(\mathcal V)$ is
\[
\xi_0(v,v') = \dfrac{d}{d \epsilon} \bigg|_{\epsilon = 0} \tilde{\bar {\mathsf m}}_{\epsilon} (v v',v').
\]
\end{proposition}

\proof We only need to show that $u_0$ is linear, but this follows from a direct computation exploiting that $R_{h_{\lambda} v} = h_\lambda \circ R_v$ for all $\lambda$, and the linearity property of $\tilde{\mathsf m}$. \endproof

We pass to $\mathsf s$-constant deformations.

\begin{proposition} [$\mathsf s$-constant deformations] Let $(\mathcal V_{\epsilon})$ be an $\mathsf s$-constant deformation of the VB-groupoid $\mathcal V$. Then
\[
\xi_0(v,v') = \dfrac{d}{d \epsilon} \bigg|_{\epsilon = 0} \tilde{\bar {\mathsf m}}_{\epsilon} (v v',v').
\]
defines a cocycle in $C^1_{\mathrm{def,lin}}(\mathcal V)$ and its cohomology class only depends on the equivalence class of the deformation.
\end{proposition}

\proof We know that $\xi_0$ is a linear cocycle. Moreover, if $\Psi_{\epsilon}: \mathcal V_{\epsilon} \to \mathcal V'_{\epsilon}$ is an equivalence of deformations, $X = \frac{d \Psi_{\epsilon}}{d \epsilon}|_{\epsilon = 0} \in C^0_{\mathrm{def,lin}}(\mathcal V)$ and we can proceed as in the proof of \cite[Lemma 5.3]{crainic:def2}. \endproof

The next step will be the proof of a linear version of Proposition \ref{prop:gen_def}. Recall from \cite{abad:ruth2} that an \emph{Ehresmann connection} on a Lie groupoid $\mathcal G \rightrightarrows M$ is a splitting of the short exact sequence (\ref{eq:ses1}) that restricts to the canonical splitting (\ref{eq:ses2}) over $M$. Notice that an Ehresmann connection on $\mathcal G$ is exactly the same as a right-horizontal lift of the VB-groupoid $T\mathcal G$. In particular, such a connection always exists.

Now let $(\mathcal V \rightrightarrows E, \mathcal G \rightrightarrows M)$ be a VB-groupoid. Then there is a diagram of morphisms of DVBs
\begin{equation}\label{eq:DVB_seq}
\begin{array}{c}
\xymatrix@C=9pt@R=12pt{& 0 \ar[rr] & & T^{\tilde {\mathsf s}} \mathcal V \ar[dl] \ar[rr] \ar[dd]|!{[dl];[dr]}{\hole} & & T \mathcal V \ar[dl] \ar[dd]|!{[dl];[dr]}{\hole} \ar[rr]^(.6){T \tilde {\mathsf s}} & & \tilde {\mathsf s}^* TE \ar[dl] \ar[dd]|!{[dl];[dr]}{\hole} \ar[rr] & & 0 \\
0 \ar[rr] & & T^{\mathsf s} \mathcal G \ar[rr] \ar[dd] & & T \mathcal G \ar[rr]^(.6){T \mathsf s} \ar[dd] & & \mathsf s^* TM \ar[rr] \ar[dd] & & 0 \\
& & & \mathcal V \ar[dl] \ar@{=}[rr]|!{[ur];[dr]}{\hole} & & \mathcal V \ar[dl] \ar@{=}[rr]|!{[ur];[dr]}{\hole} & & \mathcal V \ar[dl] & & \\
& & \mathcal G \ar@{=}[rr] & & \mathcal G \ar@{=}[rr] & & \mathcal G & &}
\end{array}
\end{equation}
where the top rows are short exact sequences of vector bundles.

\begin{definition} A \emph{linear Ehresmann connection} on $\mathcal V$ is a morphism of DVBs 
\[
(\tilde{\mathsf s}^* TE \to \mathsf s^* TM; \mathcal V \to \mathcal G) \longrightarrow (T \mathcal V \to T \mathcal G; \mathcal V \to \mathcal G)
\]
such that the maps $\tilde{\mathsf s}^* TE \to T \mathcal V$ and $\mathsf s^* TM \to T \mathcal G$ are Ehresmann connections on $\mathcal V$ and $\mathcal G$, respectively. \end{definition}

\begin{lemma} On every VB-groupoid $(\mathcal V \rightrightarrows E; \mathcal G \rightrightarrows M)$ there exists a linear Ehresmann connection. \end{lemma}

\proof First of all, choose local coordinates on $\mathcal G$ adapted to the submersion $\mathsf s: \mathcal G \to M$. Up to a translation, one can also assume that they are adapted to the immersion $\mathsf 1: M \to \mathcal G$. Using a right-decomposition $\mathcal V \cong \mathsf s^* E \oplus \mathsf t^* C$, we find fiber coordinates on $\mathcal V$ with analogous properties. Now it is easy to see that linear Ehresmann connections exist locally, and one can conclude with a partition of unity argument. \endproof

\begin{proposition} Let $\tilde{\mathcal V}$ be a deformation of $\mathcal V$. Then:
\begin{enumerate}
\item there exist transverse linear vector fields for $\tilde{\mathcal V}$;
\item if $\tilde X$ is a transverse linear vector field, then $\delta \tilde{X}$, when restricted to $\mathcal V$, induces a cocycle $\xi_0 \in C^1_{\mathrm{def,lin}}(\mathcal V)$;
\item the cohomology class of $\xi_0$ does not depend on the choice of $\tilde X$.
\end{enumerate}
\end{proposition}

\proof 
\
\begin{enumerate}
\item Take a vector field $Y$ on $\tilde M$ that projects on $\frac{d}{d\epsilon}$. Choosing a linear connection on $\tilde E \to \tilde M$, one can lift it to a linear vector field $\tilde Y$ on $\tilde E$, that obviously projects on $\frac{d}{d\epsilon}$. Now the choice of a linear Ehresmann connection on $\mathcal V$ gives a linear vector field $\tilde X$ on $\mathcal V$ that projects on $\tilde Y$, as desired.
\item Let $\tilde X$ be a transverse linear vector field. Then $\delta \tilde X \in C^1_{\mathrm{def,lin}}(\tilde{\mathcal V})$ and, by Proposition \ref{prop:gen_def}, it restricts to $\mathcal V$, so it belongs to $C^1_{\mathrm{def,lin}}(\mathcal V)$. 
\item This can be proved as in \cite[Proposition 5.12]{crainic:def2}. \endproof
\end{enumerate}

The cohomology class $[\xi_0] \in H^1_{\mathrm{def,lin}}(\mathcal V)$ is called the \emph{linear deformation class} associated to the deformation $\tilde {\mathcal V}$. From the last proposition, it follows directly that this class is also independent of the equivalence class of the deformation.

\begin{remark}\label{rem:ind_def_base}
Clearly, a deformation $\tilde{\mathcal V}$ of the VB-groupoid $\mathcal V$ induces a deformation $\tilde{\mathcal G}$ of the base groupoid $\mathcal G$. Now, it is easy to check that the projection (\ref{eq:VB_pr}) sends the linear deformation class of $\tilde{\mathcal V}$ to the deformation class of $\tilde{\mathcal G}$.
\end{remark}

Finally, we discuss the variation map associated to deformations of VB-groupoids. Let $(\tilde{\mathcal V} \rightrightarrows \tilde E; \tilde{\mathcal G} \rightrightarrows \tilde M)$ be a family of VB-groupoids over a smooth manifold $B$. Then any curve $\gamma: I \to B$ induces a deformation $\gamma^* \tilde{\mathcal V}$ of $\tilde{\mathcal V}_{\gamma(0)}$, and we have:

\begin{proposition}[The linear variation map] Let $b \in B$. For any curve $\gamma: I \to B$ with $\gamma(0) = b$, the deformation class of $\gamma^* \tilde{\mathcal V}$ at time 0 does only depend on $\dot{\gamma}(0)$. This defines a linear map
\[
\mathrm{Var}^{\tilde{\mathcal V}}_{\mathrm{lin}, b}: T_b B \to H^1_{\mathrm{def,lin}}(\tilde{\mathcal V}_b),
\]
called the \emph{linear variation map} of $\tilde{\mathcal V}$ at $b$, that makes the following diagram commutative:
\[
\xymatrix@C=45pt{
T_b B \ar[r]^-{\mathrm{Var}^{\tilde{\mathcal V}}_{\mathrm{lin}, b}} \ar[dr]_-{\mathrm{Var}^{\tilde{\mathcal G}}_b} & H^1_{\mathrm{def,lin}}(\tilde{\mathcal V}_b) \ar[d] \\
& H^1_{\mathrm{def}}(\tilde{\mathcal G}_b)
}
\]
\end{proposition}

\proof The first statement is proved as in \cite[Proposition 5.15]{crainic:def2}, while the second statement trivially follows from Remark \ref{rem:ind_def_base}. \endproof

\subsubsection{Deformations of the dual VB-groupoid}

We conclude this section noticing that the linear deformation cohomology of a VB-groupoid is canonically isomorphic to that of its dual.

\begin{theorem} Let $(\mathcal V \rightrightarrows E; \mathcal G \rightrightarrows M)$ be a VB-groupoid. Then there is a canonical isomorphism
\begin{equation}\label{eq:dual_def}
H_{\mathrm{def,lin}}(\mathcal V) \cong H_{\mathrm{def,lin}}(\mathcal V^*)
\end{equation}
of $H(\mathcal G)$-modules.
\end{theorem}

\proof Using Proposition \ref{prop:cotangent_VB} and Lemma \ref{prop:proj,lin;lin,lin}, we get:
\[
H_{\mathrm{def,lin}}(\mathcal V)\cong H_{\mathrm{proj,lin}}(T^* \mathcal V)[1]  \cong H_{\mathrm{lin,lin}}(T^* \mathcal V)[1].
\]
For the same reason, $H_{\mathrm{def,lin}}(\mathcal V^*) \cong H_{\mathrm{lin,lin}}(T^* \mathcal V^*)[1]$. But we have already noticed that $T^* \mathcal V \cong T^* \mathcal V^*$ as double vector bundles and Lie groupoids, so we obtain (\ref{eq:dual_def}). \endproof

\subsection{The linearization map}\label{sec:linearization_1}

Let $(\mathcal V \rightrightarrows E; \mathcal G \rightrightarrows M)$ be a VB-groupoid. We have shown that deformations of the VB-groupoid structure are controlled by a subcomplex $C_{\mathrm{def}, \mathrm{lin}} (\mathcal V)$ of the deformation complex $C_{\mathrm{def}}(\mathcal V)$ of the top Lie groupoid $\mathcal V \rightrightarrows E$. In this section, we prove that there is a canonical splitting of the inclusion $C_{\mathrm{def,lin}}(\mathcal V) \hookrightarrow C_{\mathrm{def}}(\mathcal V)$ in the category of cochain complexes, the \emph{linearization map}. This will imply, in particular, that the inclusion induces an injection in cohomology $H_{\mathrm{def,lin}}(\mathcal V) \hookrightarrow H_{\mathrm{def}}(\mathcal V)$.

The procedure we are going to describe is analogous to the one we used in \cite{lapastina:def} to define the linearization of sections of a DVB. Following Remark \ref{rmk:hom_str}, denote by $h$ the homogeneity structure of $\mathcal V \to \mathcal G$. For every $\lambda > 0$, $h_\lambda$ is a groupoid automorphism of $\mathcal V \rightrightarrows E$. If $k = -1$, by definition the action of $h_\lambda$ on $C^{-1}_{\mathrm{def}}(\mathcal V) = \Gamma(W,E)$ coincides with that induced by  the homogeneity structure of $W \to A$. If $k \geq 0$, $C^k_{\mathrm{def}}(\mathcal V) \subset \Gamma(\tilde p_{k+1}^* T \mathcal V, \mathcal V^{(k+1)})$, and a direct computation shows that the action of $h_\lambda$ on $C^k_{\mathrm{def}}(\mathcal V)$ coincides with that induced by the homogeneity structure of the vector bundle $\tilde p_{k+1}^* T \mathcal V \to p_{k+1}^* T \mathcal G$ on sections of $p_{k+1}^* T \mathcal V \to \mathcal V^{(k+1)}$. Therefore, by \cite[Propositions 1.4.1, 1.4.2]{lapastina:def}, the limits
\[
\begin{aligned}
\tilde c_{\mathrm{core}} := & \lim_{\lambda \to 0} \big( \lambda \cdot h_\lambda^* \tilde c \big) \\
\tilde c_{\mathrm{lin}} := & \lim_{\lambda \to 0} \big( h_\lambda^* \tilde c - \lambda^{-1} \cdot \tilde c_{\mathrm{core}} \big) \\
\end{aligned}
\]
are well-defined. The latter equation defines a linear map
\begin{equation}\label{eq:linearization}
\mathrm{lin}: C_{\mathrm{def}}(\mathcal V) \to C_{\mathrm{def,lin}}(\mathcal V), \quad \tilde c \mapsto \tilde c_{\mathrm{lin}},
\end{equation}
which we call the \emph{linearization map}.

\begin{theorem}\label{prop:linearization}
The map (\ref{eq:linearization}) is a cochain map that splits the inclusion $C_{\mathrm{def,lin}}(\mathcal V) \hookrightarrow C_{\mathrm{def}}(\mathcal V)$. Hence $C_{\mathrm{def,lin}}(\mathcal V)$ is a direct summand of $C_{\mathrm{def}}(\mathcal V)$ and the inclusion induces an injection in cohomology:
\[
H_{\mathrm{def,lin}}(\mathcal V) \hookrightarrow H_{\mathrm{def}}(\mathcal V).
\]
\end{theorem}

\proof We only need to prove that the linearization map respects the differential. As $h_\lambda$ is an automorphism of $\mathcal V \rightrightarrows E$ for every $\lambda$, we have that $h_\lambda^*$ commutes with $\delta$. It is also clear that $\delta$ preserves limits, so we compute
\[
\begin{aligned}
(\delta \tilde c)_{\mathrm{core}} = & \lim_{\lambda \to 0} \big( \lambda \cdot h_\lambda^* (\delta \tilde c) \big) = \lim_{\lambda \to 0} \big( \lambda \cdot \delta (h_\lambda^* \tilde c) \big) = \delta \big( \lim_{\lambda \to 0} \lambda \cdot \delta (h_\lambda^* \tilde c) \big) = \delta \tilde c_{\mathrm{core}}, \\
(\delta \tilde c)_{\mathrm{lin}} = & \lim_{\lambda \to 0} \big( h_\lambda^* (\delta \tilde c) - \lambda^{-1} (\delta \tilde c)_{\mathrm{core}} \big) = \delta \bigg( \lim_{\lambda \to 0} (h_\lambda^* \tilde c - \lambda^{-1} \delta \tilde c_{\mathrm{core}}) \bigg) = \delta \tilde c_{\mathrm{lin}}
\end{aligned}
\]
and we are done. \endproof 

\begin{remark}\label{rmk:cotangent} Applying the isomorphisms (\ref{eq:iso}) and (\ref{eq:iso_lin}), we obtain that \emph{$H_{\mathrm{proj,lin}}(T^* \mathcal V)$ is a direct summand of $H_{\mathrm{proj},\bullet}(T^* \mathcal V)$}: it identifies with classes in $H_{\mathrm{proj},\bullet}(T^* \mathcal V)$ which can be represented by cochains that are linear over $\mathcal V^*$. \end{remark}

Finally, we discuss a first consequence of Theorem \ref{prop:linearization}. We call $\widehat C_{\mathrm{def,lin}}(\mathcal V) := C_{\mathrm{def,lin}}(\mathcal V) \cap \widehat C_{\mathrm{def}}(\mathcal V)$ the \emph{linear normalized deformation complex} of $\mathcal V$.

\begin{proposition} The inclusion $\widehat C_{\mathrm{def,lin}}(\mathcal V) \hookrightarrow C_{\mathrm{def,lin}}(\mathcal V)$ is a quasi-isomorphism. \end{proposition}

\proof Take $c \in C^k_{\mathrm{def,lin}}(\mathcal V)$. By Proposition \ref{prop:norm}, there exist $\widehat c \in \widehat C^k_{\mathrm{def}}(\mathcal V)$ and $c' \in C^{k-1}_{\mathrm{def}}(\mathcal V)$ such that $c - \widehat c = \delta c'$. Applying the linearization map, we get $c - \widehat c_{\mathrm{lin}} = \delta c'_{\mathrm{lin}}$ and $\widehat c_{\mathrm{lin}} \in \widehat C^k_{\mathrm{def,lin}}(\mathcal V)$, as desired. \endproof

Other applications of the linearization map will be considered in the next sections.

\subsection{The van Est map} \label{sec:van_est}

The van Est theorem is a classical result relating the differentiable cohomology of a Lie group and the Chevalley-Eilenberg cohomology of its Lie algebra \cite{est:group_coh, est:alg_coh}. It was later extended to differentiable cohomology \cite{weinstein:extensions, crainic:coh} and deformation cohomology \cite{crainic:def2} of a Lie groupoid, and to the VB-cohomology of a VB-groupoid \cite{cabrera:hom}. In this subsection, we want to prove an analogous theorem for the linear deformation cohomology of a VB-groupoid.

Let $\mathcal G \rightrightarrows M$ be a Lie groupoid, and let $A \Rightarrow M$ be its Lie algebroid. The normalized deformation complex of $\mathcal G$ and the deformation complex of $A$ are intertwined by the \emph{van Est map}, defined as follows. Given a section $\alpha \in \Gamma(A)$, we define a map $R_\alpha: \widehat C^k_{\mathrm{def}}(\mathcal G) \to \widehat C^{k-1}_{\mathrm{def}}(\mathcal G)$ by
\begin{equation}\label{eq:R_alpha_1}
R_\alpha(c) = [c, \overrightarrow{\alpha}]|_M
\end{equation}
if $k = 0$, and
\begin{equation}\label{eq:R_alpha_2}
R_\alpha(c) (g_1, \dots, g_k) = (-1)^{k} \dfrac{d}{d \epsilon} \bigg|_{\epsilon = 0} c(g_1, \dots, g_k, \Phi^\alpha_\epsilon(\mathsf s (g_k))^{-1})
\end{equation}
if $k > 0$, where $\Phi^\alpha_\epsilon (x) = \Phi^{\overrightarrow{\alpha}}_\epsilon (\mathsf 1_x)$ for every $x \in M$: the image of $1_x$ under the flow $\{\Phi^{\overrightarrow{\alpha}}_\epsilon\}$ of the right invariant vector field $\overrightarrow{\alpha}$ associated to $\alpha$. Then the van Est map 
\begin{equation}\label{eq:VE}
\mathrm{VE}: \widehat C_{\mathrm{def}}(\mathcal G) \to C_{\mathrm{def}}(A)
\end{equation}
is given by:
\begin{equation}\label{eq:VE_formula}
\mathrm{VE}(c)(\alpha_0, \dots, \alpha_k) = \sum_{\tau \in S_{k+1}} (-1)^\tau (R_{\alpha_{\tau(k)}} \circ \dots \circ R_{\alpha_{\tau(0)}}) (c).
\end{equation}

\begin{theorem}[\cite{crainic:def2}]\label{prop:van_est}  The van Est map $\mathrm{VE}$ is a cochain map. Moreover, if $\mathcal G$ has $k$-connected $\mathsf s$-fibers, it induces an isomorphism in cohomology in all degrees $p < k$. \end{theorem}

We are going to prove an analogous theorem for the linear deformation complex of a VB-groupoid. To do this, we need a simple preliminary lemma.

Let $\mathcal G \rightrightarrows M$ be a Lie groupoid, and let $A \Rightarrow M$ be its Lie algebroid.  We know that the group $\mathrm{Aut}(\mathcal G)$ of automorphisms of $\mathcal G$ acts on $\widehat C_{\mathrm{def}}(\mathcal G)$. Clearly, it also acts on $C_{\mathrm{def}}(A)$, by
\[
\Psi^* c := \psi^* c.
\]
As expected, we have the following
\begin{lemma}\label{lem:VE_equiv}
 The van Est map (\ref{eq:VE}) is equivariant with respect to the action of $\mathrm{Aut}(\mathcal G)$. \end{lemma}

\proof We will prove that
\begin{equation}\label{eq:R_equiv}
\Psi^*(R_\alpha (c)) = R_{\psi^* \alpha}(\Psi^* c)
\end{equation}
for all $c \in \widehat C^k_{\mathrm{def}}(\mathcal G), \alpha \in \Gamma(A), \Psi \in \mathrm{Aut}(\mathcal G)$. If $k = 0$, we have
\[
\Psi^*(R_\alpha(c)) = \Psi^*([c, \overrightarrow \alpha]|_M) = (\Psi^*[c, \overrightarrow \alpha])|_M = [\Psi^* c, \Psi^* \overrightarrow \alpha]|_M = [\Psi^*c, \overrightarrow{\psi^* \alpha}]|_M = R_{\psi^* \alpha}(\Psi^* c).
\]
Now, let $k > 0$. Then
\begin{equation}\label{eq:Psi^*R}
\begin{aligned}
\Psi^*(R_\alpha(c)) (g_1, \dots, g_k) = & T \Psi^{-1} (R_\alpha (c) (\Psi(g_1), \dots, \Psi(g_k))) \\ 
= & T\Psi^{-1} \bigg( (-1)^{k} \dfrac{d}{d \epsilon} \bigg|_{\epsilon = 0} c(\Psi(g_1), \dots, \Psi(g_k), \Phi^\alpha_\epsilon(\mathsf s (\Psi(g_k)))^{-1}) \bigg)
\end{aligned}
\end{equation} 
Using (\ref{eq:inv_vf}), we compute:
\[
\begin{aligned}
\Phi^\alpha_\epsilon(\mathsf s(\Psi(g_k))) & = \Phi^{\overrightarrow{\alpha}}_\epsilon(\mathsf 1_{\mathsf s(\Psi(g_k))}) = \Phi^{\overrightarrow{\alpha}}_\epsilon(\Psi(\mathsf 1_{\mathsf s(g_k)})) \\ & = \Psi(\Phi^{\Psi^*\overrightarrow{\alpha}}_\epsilon(\mathsf 1_{\mathsf s(g_k)})) = \Psi(\Phi_\epsilon^{\overrightarrow{\psi^* \alpha}}(\mathsf 1_{\mathsf s(g_k)})) = \Psi(\Phi^{\psi^* \alpha}_\epsilon(\mathsf s(g_k))).
\end{aligned}
\]
So, from (\ref{eq:Psi^*R}) we have:
\[
\begin{aligned}
\Psi^*(R_\alpha(c)) (g_1, \dots, g_k) = & (-1)^{k}  \dfrac{d}{d \epsilon} \bigg|_{\epsilon = 0} \Psi^{-1}(c(\Psi(g_1), \dots, \Psi(g_k), \Psi(\Phi_\epsilon^{\psi^* \alpha} (\mathsf s(g_k)))^{-1})) \\
= & (-1)^{k} \dfrac{d}{d\epsilon} \bigg|_{\epsilon = 0} (\Psi^* c)(g_1, \dots, g_k, \Phi_\epsilon^{\psi^* \alpha}(\mathsf s(g_k))^{-1}) \\
= & R_{\psi^* \alpha}(\Psi^* c)(g_1, \dots, g_k).
\end{aligned}
\]
Finally, by applying repeatedly formula (\ref{eq:R_equiv}) in (\ref{eq:VE_formula}), we obtain
\[
\Psi^*(\mathrm{VE}(c)) = \mathrm{VE}(\Psi^* c)
\]
for every $c \in \widehat C_{\mathrm{def}}(\mathcal G)$, as desired. \endproof

Now we are ready for the main theorem of this section. Notice that the first part of the following statement has been already proved in \cite{etv:infinitesimal} in the special case where the rank of $E$ is greater than zero. Here we provide an alternative proof exploiting Lemma \ref{lem:VE_equiv} which is valid in all cases.

\begin{theorem}[Linear van Est map] \label{prop:lin_van_est} Let $(\mathcal V \rightrightarrows E; \mathcal G \rightrightarrows M)$ be a VB-groupoid, and let $(W \Rightarrow E; A \Rightarrow M)$ be its VB-algebroid. Then the van Est map for the Lie groupoid $\mathcal V \rightrightarrows E$ restricts to a cochain map
\[
\mathrm{VE}: \widehat C_{\mathrm{def,lin}}(\mathcal V) \to C_{\mathrm{def,lin}}(W),
\]
which we call the \emph{linear van Est map}. If $\mathcal G$ has $k$-connected $\mathsf s$-fibers, this map induces an isomorphism in cohomology in all degrees $p < k$. \end{theorem}

\proof As before, we denote by $h$ the homogeneity structure of $\mathcal V \to \mathcal G$. Then the last proposition shows that
\[
h_\lambda^* (\mathrm{VE}(\tilde c)) = \mathrm{VE}(h_\lambda^* \tilde c)
\]
for every $\tilde c \in \widehat C_{\mathrm{def}}(\mathcal V)$, $\lambda > 0$. But the VB-algebroid automorphism corresponding to $h_\lambda$ is exactly the one induced by the homogeneity structure of $W \to A$, so the last equation implies that the van Est map preserves linear cochains.

Now we prove the second part of the theorem. First, we have to observe that the van Est map commutes with linearization, i.e.
\begin{equation}\label{eq:VE_lin}
\mathrm{VE}(\tilde c)_{\mathrm{lin}} = \mathrm{VE}(\tilde c_{\mathrm{lin}}).
\end{equation}
To see this, take $\tilde c \in C^k_{\mathrm{def,lin}}(\mathcal V)$ and $w_0, \dots, w_k \in \Gamma(W,E)$ and compute
\[
\begin{aligned}
\mathrm{VE}(\tilde c)_{\mathrm{lin}} (w_0, \dots, w_k) & = \lim_{\lambda \to 0} \big( h_\lambda^* (\mathrm{VE}(\tilde c)) - \lambda^{-1} \mathrm{VE}(\tilde c_{\mathrm{core}}) \big) (w_0, \dots, w_k) \\
& = \lim_{\lambda \to 0} \mathrm{VE}(h_\lambda^* \tilde c - \lambda^{-1} \tilde c_{\mathrm{core}}) (w_0, \dots, w_k) \\
& = \lim_{\lambda \to 0} \sum_{\tau \in S_{k+1}} (-1)^\tau (R_{w_{\tau(k)}} \circ \dots \circ R_{w_{\tau(0)}}) (h_\lambda^* \tilde c - \lambda^{-1} \tilde c_{\mathrm{core}}).
\end{aligned}
\]
We would like to swap the limit and the derivatives that appear in definitions (\ref{eq:R_alpha_1}) and (\ref{eq:R_alpha_2}). This is ultimately possible because of smoothness,
and we get
\[
\begin{aligned}
\mathrm{VE}(\tilde c)_{\mathrm{lin}}(w_0, \dots, w_k) & = \sum_{\tau \in S_{k+1}} (-1)^\tau (R_{w_{\tau(k)}} \circ \dots \circ R_{w_{\tau(0)}}) \big( \lim_{\lambda \to 0} (h_\lambda^* \tilde c - \lambda^{-1} \tilde c_{\mathrm{core}}) \big) \\
& = \sum_{\tau \in S_{k+1}} (-1)^\tau (R_{w_{\tau(k)}} \circ \dots \circ R_{w_{\tau(0)}}) (\tilde c_{\mathrm{lin}}) \\
& = \mathrm{VE}(\tilde c_{\mathrm{lin}}) (w_0, \dots, w_k)
\end{aligned}
\]
as desired.

Finally, suppose that $\mathcal G$ has $k$-connected $\mathsf s$-fibers. Then $\mathcal V$ has $k$-connected $\tilde{\mathsf s}$-fibers (they are vector bundles over the $\mathsf s$-fibers of $\mathcal G$). Take $p < k$. We want to prove that the induced map
\[
\mathrm{VE}: H^p_{\mathrm{def,lin}}(\mathcal V) \to H^p_{\mathrm{def,lin}}(W)
\]
is an isomorphism. If $\tilde c \in C^p_{\mathrm{def,lin}}(\mathcal V)$ is closed, we denote $[\tilde c]$ its class in $H^p_{\mathrm{def}}(\mathcal V)$ and $[\tilde c]_{\mathrm{lin}}$ its class in $H^p_{\mathrm{def,lin}}(\mathcal V)$; we use an analogous notation for $W$.

First, suppose that $[\mathrm{VE}(\tilde c)]_{\mathrm{lin}} = 0$. Then $[\mathrm{VE}(\tilde c)] = 0$ and, by Theorem \ref{prop:van_est}, $[\tilde c] = 0$, i.e. $\tilde c = \delta \tilde \gamma$ for some $\tilde \gamma \in C^{p-1}_{\mathrm{def}}(\mathcal V)$. Applying the linearization map, we get $\tilde c = \delta \tilde \gamma_{\mathrm{lin}}$ and $\tilde \gamma_{\mathrm{lin}} \in \widehat C^{p-1}_{\mathrm{def,lin}}(\mathcal V)$, i.e.~$\mathrm{VE}$ is injective in degree $p$ cohomology. To conclude, take $[c]_{\mathrm{lin}} \in H^p_{\mathrm{def,lin}}(W)$. Then $[c] \in H^p_{\mathrm{def}}(W)$, so $[c] = [\mathrm{VE}(\tilde c)]$, i.e. $c - \mathrm{VE}(\tilde c) = \delta \gamma$ for some $\tilde c \in \widehat C^p_{\mathrm{def}}(\mathcal V), \gamma \in C^{p-1}_{\mathrm{def}}(W)$. Applying again the linearization map and using (\ref{eq:VE_lin}), we get $c - \mathrm{VE}(\tilde c_{\mathrm{lin}}) = \delta \gamma_{\mathrm{lin}}$, i.e.~$\mathrm{VE}$ is also surjective in degree $p$ cohomology, as desired. \endproof

\subsection{Morita invariance} \label{sec:morita}

The notion of Morita equivalence of VB-groupoids first appears in \cite{delhoyo:morita}. In that reference, the authors prove that the VB-cohomologies of Morita equivalent VB-groupoids are isomorphic. As a corollary, they give a conceptual and very simple proof of the fact, first appeared in \cite{crainic:def2}, that Morita equivalent Lie groupoids have isomorphic deformation cohomologies. This second result means that the deformation cohomology of a Lie groupoid is in fact an invariant of the associated differentiable stack.

In this paragraph, we want to prove an analogous result for the linear deformation cohomology of a VB-groupoid. We start recalling the necessary definitions. Let $\mathcal G_1\rightrightarrows M_1$ and $\mathcal G_2 \rightrightarrows M_2$ be Lie groupoids. 

\begin{definition}\label{def:Morita_map} A morphism of Lie groupoids $\Psi: \mathcal G_1 \to \mathcal G_2$ over a smooth map $F: M_1 \to M_2$ is a \emph{Morita map} (or a \emph{weak equivalence}) if it is
\begin{enumerate}
\item \emph{fully faithful}, i.e. the diagram
\[
\begin{array}{c}
\xymatrix@C=45pt{\mathcal G_1 \ar[r]^{\Psi} \ar[d]_{(\mathsf s_1, \mathsf t_1)} & \mathcal G_2 \ar[d]^{(\mathsf s_2, \mathsf t_2)} \\
M_1 \times M_1 \ar[r]_{F \times F} & M_2 \times M_2}
\end{array}
\]
is a pull-back diagram, and
\item \emph{essentially surjective}, i.e. the map $\mathcal G_2 \tensor*[_{\mathsf s_1}]{\times}{_F}  M_1 \to M_2, (g,x) \mapsto \mathsf t_2(g)$ is a surjective submersion.
\end{enumerate}

Two groupoids $\mathcal G_1$ and $\mathcal G_2$ are said to be \emph{Morita} (or \emph{weakly}) \emph{equivalent} if and only if there exist a Lie groupoid $\mathcal H$ and Morita maps $\Psi_1: \mathcal H \to \mathcal G_1$, $\Psi_2: \mathcal H \to \mathcal G_2$.
\end{definition}

Notice that Definition \ref{def:Morita_map}(1) of fully faithful morphism is slightly different from the one in \cite{delhoyo:morita}, where an additional property is required. However, it is easy to see that for an essentially surjective morphism the two definitions are equivalent. In the literature, Morita equivalence is often expressed in terms of principal bibundles: we refer to \cite{mrcun:stability} for this notion and many more details.

Now, let $(\mathcal V_1 \rightrightarrows E_1; \mathcal G_1 \rightrightarrows M_1)$ and $(\mathcal V_2 \rightrightarrows E_2; \mathcal G_2 \rightrightarrows M_2)$ be VB-groupoids. A VB-groupoid morphism $\Psi: \mathcal V_1 \to \mathcal V_2$ is a \emph{VB-Morita map} \cite{delhoyo:morita} if the Lie groupoid morphism $\Psi$ is a Morita map. The VB-groupoids $\mathcal V_1$ and $\mathcal V_2$ are \emph{Morita equivalent} if there exist a VB-groupoid $\mathcal W$ and VB-Morita maps $\mathcal W \to \mathcal V_1$, $\mathcal W \to \mathcal V_2$.

Here are some basic properties of VB-Morita maps.

\begin{proposition}[{\cite[Corollary 3.7]{delhoyo:morita}}]\label{prop:VBM_1} Let $\Psi: \mathcal G_1 \to \mathcal G_2$ be a Morita map and let $\mathcal V$ be a VB-groupoid over $\mathcal G_2$. Then the canonical map $\Psi^* \mathcal V \to \mathcal V$ is VB-Morita. \end{proposition}

\begin{proposition}[{\cite[Corollary 3.8]{delhoyo:morita}}]\label{prop:VBM_2} Let $\Psi: \mathcal G_1 \to \mathcal G_2$ be a Morita map. Then its tangent map $T \Psi: T \mathcal G_1 \to T \mathcal G_2$ is a VB-Morita map. \end{proposition} 

\begin{proposition}[{\cite[Corollary 3.9]{delhoyo:morita}}]\label{prop:VBM_3} A map $\Psi: \mathcal V_1 \to \mathcal V_2$ over the identity is VB-Morita if and only if its dual is so. \end{proposition}

Morita invariance of the VB-cohomology is expressed by the following theorem.

\begin{theorem} [{\cite[Theorem 4.2]{delhoyo:morita}}] Let $\Psi: \mathcal V_1 \to \mathcal V_2$ be a VB-Morita map. Then $\Psi^*: H_{\mathrm{proj}}(\mathcal V_2) \to H_{\mathrm{proj}}(\mathcal V_1)$ is an isomorphism. \end{theorem}

Now we are ready to prove Morita invariance of the linear deformation cohomology.

\begin{theorem} Let $\Psi: \mathcal V_1 \to \mathcal V_2$ be a VB-Morita map. Then $H_{\mathrm{def,lin}}(\mathcal V_1) \cong H_{\mathrm{def,lin}}(\mathcal V_2)$. \end{theorem}

\proof It is enough to show that $H_{\mathrm{proj,lin}}(T^* \mathcal V_1) \cong H_{\mathrm{proj,lin}}(T^* \mathcal V_2)$. To do this we will use Propositions \ref{prop:VBM_1}--\ref{prop:VBM_3} and linearization.

Recall that both $T^* \mathcal V_1$ and $T^* \mathcal V_2$ have two VB-groupoid structures, as discussed in Subsection \ref{sec:lin_def}, and observe that $\Psi^* (T^* \mathcal V_2)$ possesses also two VB-groupoid structures, that fit in the following commuting diagram:
\[
\begin{array}{c}
\xymatrix@C=9pt@R=12pt{
& \Psi^* (T^* \mathcal V_2) \ar[dl] \ar@<0.4ex>[rr] \ar@<-0.4ex>[rr] \ar[dd]|!{[dl];[dr]}{\hole} & & \Psi^* (W_2^* |_{E_2}) \ar[dl] \ar[dd] \\
\mathcal V_1 \ar@<0.4ex>[rr] \ar@<-0.4ex>[rr] \ar[dd] & & E_1 \ar[dd]\\
& \Psi^* \mathcal V_2^* \ar[dl] \ar@<0.4ex>[rr]|!{[ur];[dr]}{\hole} \ar@<-0.4ex>[rr]|!{[ur];[dr]}{\hole} & & \Psi^* C_2^* \ar[dl] \\
\mathcal G_1 \ar@<0.4ex>[rr] \ar@<-0.4ex>[rr] & & M_1 
}
\end{array}.
\]

We denote by $C_{\mathrm{proj},\bullet}(\Psi^* (T^* \mathcal V_2))$ the VB-complex of the VB-groupoid upstairs and by $C_{\mathrm{proj,lin}}(\Psi^* (T^* \mathcal V_2))$ its subcomplex of cochains that are linear with respect to the vertical projections. As usual, we denote their cohomologies by $H_{\mathrm{proj},\bullet}(\Psi^* (T^* \mathcal V_2))$ and $H_{\mathrm{proj,lin}}(\Psi^* (T^* \mathcal V_2))$ respectively. In analogy to Theorem \ref{prop:linearization} and Remark \ref{rmk:cotangent}, one can prove that there is a linearization map
\[
\mathrm{lin}: C_{\mathrm{proj},\bullet}(\Psi^*(T^* \mathcal V_2)) \to C_{\mathrm{proj,lin}}(\Psi^* (T^* \mathcal V_2))
\]
that splits the inclusion $C_{\mathrm{proj,lin}}(\Psi^* (T^* \mathcal V_2)) \hookrightarrow C_{\mathrm{proj},\bullet}(\Psi^*(T^* \mathcal V_2))$ in the category of cochain complexes. Hence $H_{\mathrm{proj,lin}}(\Psi^* (T^* \mathcal V_2))$ embeds in $H_{\mathrm{proj},\bullet}(\Psi^* (T^* \mathcal V_2))$ as a direct summand.

The map $\Psi$ is VB-Morita, so, from Proposition \ref{prop:VBM_2}, $T \Psi: T \mathcal V_1 \to T \mathcal V_2$ is VB-Morita again. It follows from Propositions \ref{prop:VBM_1} and \ref{prop:VBM_3} that the dual map $(T \Psi)^*: \Psi^* (T^* \mathcal V_2) \to T^* \mathcal V_1$ is VB-Morita as well.  Remember that also the canonical map $\Psi^* (T^* \mathcal V_2) \to T^* \mathcal V_2$ is VB-Morita (Proposition \ref{prop:VBM_1} again). As a result, we get isomorphisms in VB-cohomology:
\begin{equation}\label{eq:proj_coh}
H_{\mathrm{proj},\bullet}(T^* \mathcal V_1) \overset{\cong}{\longrightarrow} H_{\mathrm{proj},\bullet}(\Psi^*(T^* \mathcal V_2)) \overset{\cong}{\longleftarrow} H_{\mathrm{proj},\bullet}(T^* \mathcal V_2).
\end{equation}
Now, the maps 
\[
T^* \mathcal V_1 \overset{(T\Psi)^*}{\longleftarrow} \Psi^* (T^* \mathcal V_2) \longrightarrow T^* \mathcal V_2
\]
are also DVB morphisms. This implies, on one hand, that the maps (\ref{eq:proj_coh}) preserve linear cohomologies, on the other hand that they commute with the respective linearization maps. From this last property we deduce, as in Theorem \ref{prop:lin_van_est}, that the maps (\ref{eq:proj_coh}) induce isomorphisms on linear cohomologies, so
\[
H_{\mathrm{proj,lin}}(T^* \mathcal V_1) \cong H_{\mathrm{proj,lin}}(\Psi^* (T^* \mathcal V_2)) \cong H_{\mathrm{proj,lin}}(T^* \mathcal V_2)
\]
as desired. \endproof

\section{Examples and applications}\label{Sec:examples}

In this section we provide several examples. Examples in Subsections \ref{sec:VB_grp}, \ref{sec:fol_grpd} and \ref{sec:Lie_vect} parallel the analogous examples in \cite{crainic:def2}, connecting our linear deformation cohomology to known cohomologies, while examples in Subsections \ref{sec:2-vect} and \ref{sec:tangent-VB} are specific to VB-groupoids. The infinitesimal counterparts of all these examples were discussed in our previous paper \cite{lapastina:def}.

\subsection{VB-groups and their duals}\label{sec:VB_grp}

A \emph{VB-group} is a \emph{vector bundle object in the category of Lie groups}. In other words, it is a VB-groupoid of the form
\[
\begin{array}{r}
\xymatrix{
H \ar[d] \ar@<0.4ex>[r] \ar@<-0.4ex>[r] & 0 \ar[d] \\
G \ar@<0.4ex>[r] \ar@<-0.4ex>[r] & \ast}
\end{array}.
\]
In particular, $H$ and $G$ are Lie groups. Let $C := \ker (H \to G)$ be the core of $(H \rightrightarrows 0; G \rightrightarrows \ast)$. It easily follows from the definition of VB-groupoid that
\begin{itemize}
\item $C$ is a representation of $G$,
\item $H \cong G \ltimes C$ is the semidirect product Lie group,
\item $H \cong G \ltimes C \to G$ is the projection onto the first factor.
\end{itemize}
It is then natural to study the relationship between the linear deformation complex of $H$ and the classical complex $C(G, \operatorname{End} C)$ (of the Lie group $G$ with coefficients in the representation $\operatorname{End} C$)  that controls deformations of the $G$-module $C$ \cite{nijenhuis:def}. To do this, we notice that the dual of $H$ is the VB-groupoid
\[
\begin{array}{r}
\xymatrix{
G \ltimes C^* \ar[d] \ar@<0.4ex>[r] \ar@<-0.4ex>[r] & C^* \ar[d] \\
G \ar@<0.4ex>[r] \ar@<-0.4ex>[r] & \ast}
\end{array},
\]
i.e.~it is the action VB-groupoid associated to the dual representation of $G$ on $C^*$. In particular, it is a trivial-core VB-groupoid, so there is a short exact sequence of cochain complexes:
\[
0 \longrightarrow C(G, \operatorname{End} C^*) \longrightarrow C_{\mathrm{def,lin}}(G \ltimes C^*) \longrightarrow C_{\mathrm{def}}(G) \longrightarrow 0.
\]
But $C(G, \operatorname{End} C^*) \cong C(G, \operatorname{End} C)$ canonically, so the latter is recovered as the subcomplex of $C_{\mathrm{def,lin}}(G \ltimes C^*)$ controlling deformations of the representation $C^*$ that fix the Lie group structure on $G$. Moreover, $H_{\mathrm{def,lin}}(G \ltimes C^*) \cong H_{\mathrm{def,lin}}(H)$, so there is a long exact sequence in cohomology:
\[
\cdots \longrightarrow H^k(G, \operatorname{End} C) \longrightarrow H^k_{\mathrm{def}, \mathrm{lin}}(H) \longrightarrow H^k_{\mathrm{def}}(G) \longrightarrow H^{k+1}(G, \operatorname{End} C) \longrightarrow \cdots .
\]

\subsection{2-vector spaces}\label{sec:2-vect}

A \emph{2-vector space} is a \emph{(Lie) groupoid object in the category of vector spaces}. In other words, it is a VB-groupoid of the form
\[
\begin{array}{r}
\xymatrix{
V_1 \ar[d] \ar@<0.4ex>[r] \ar@<-0.4ex>[r] & V_0 \ar[d] \\
\ast \ar@<0.4ex>[r] \ar@<-0.4ex>[r] & \ast}
\end{array}.
\]

In \cite{baez:higher} it is proved that, if $\mathsf s$ and $\mathsf t$ are the source and the target maps of $V_1 \rightrightarrows V_0$, $C = \ker \mathsf s$, $\partial = \mathsf t|_C: C \to V_0$, then $V_1 \rightrightarrows V_0$ is canonically isomorphic to the action groupoid $C \ltimes V_0 \rightrightarrows V_0$, where $C$ acts on $V_0$ by
\[
c \cdot v = \partial c + v.
\] 
Notice that $C$ does not act by linear isomorphisms, but by translations. We will identify $V_1$ with $C \ltimes V_0$.

Now we compute the linear deformation complex of $(V_1 \rightrightarrows V_0; * \rightrightarrows *)$. First of all, as $V_1 \rightrightarrows V_0$ is an action groupoid, Equation (\ref{eq:nerve_act}) yields an isomorphism $(C \ltimes V_0)^{(k)} \cong C^k \oplus V_0$. We will understand this isomorphism.  Recall also that $\mathsf s: C \ltimes V_0 \to V_0$ is just the projection onto the second factor. It follows that, for $k \geq 0$, $C^k_{\mathrm{def,lin}}(V_1)$ is the set of linear maps $C^{k+1} \oplus V_0 \to C \oplus V_0$ such that the second component does not depend on the first arrow $(c_1, (c_2 + \dots + c_k) \cdot v)$, hence these maps are equivalent to couples of linear maps $C^{k+1} \oplus V_0 \to C$ and $C^{k} \oplus V_0 \to V_0$. So
\begin{equation}\label{eq:2-vect}
C^k_{\mathrm{def,lin}}(V_1) \cong \mathrm{Hom}(C^{k+1} \oplus V_0, C) \oplus \mathrm{Hom}(C^k \oplus V_0, V_0).
\end{equation}
We will identify a deformation cochain $\gamma$ with the corresponding pair of linear maps $(\gamma_1, \gamma_2)$. A direct computation shows that the differential is given by the following formulas. For $\gamma \in C^{-1}_{\mathrm{def,lin}}(V_1)$
\[
\delta \gamma = (\gamma \circ \partial, \partial \circ \gamma),
\]
and, for $\gamma \in C^k_{\mathrm{def,lin}}(V_1)$, $k \geq 0$,
\[
\delta \gamma = (\Gamma_1, \Gamma_2),
\]
where
\[
\begin{aligned}
& \Gamma_1 (c_0, \ldots, c_{k+1}, v) \\
 = & -\gamma_1 (c_0, 0, \dots, 0) + \sum_{i=1}^k (-1)^{i-1} \gamma_1 (c_0, \dots, c_i + c_{i+1}, \dots, c_{k+1}, v)  \\
 &+ (-1)^k \gamma_1 (c_0, \dots, c_k, c_{k+1} \cdot v),
\end{aligned}
\]
and 
\[
\begin{aligned}
& \Gamma_2 (c_1, \ldots, c_{k+1}, v) \\
= & - \gamma_1 (c_1, \dots, c_{k+1}, v) \cdot \gamma_2 (c_2, \dots, c_{k+1}, v) \\
& + \sum_{i=1}^k (-1)^{i-1} \gamma_2 (c_1, \dots, c_i + c_{i+1}, \dots, c_{k+1}, v) + (-1)^k \gamma_2 (c_1, \dots, c_k, c_{k+1} \cdot v).
\end{aligned}
\]

We further notice that, if $\gamma \in C^k_{\mathrm{def,lin}}(V_1)$, $k \geq 0$, then $\gamma = (\gamma_1, \gamma_2)$ belongs to the normalized deformations subcomplex $\widehat C_{\mathrm{def,lin}}(V_1)$ if and only if
\begin{center}
\begin{tabular}{cc}
$\gamma_1 (c_0, \dots, \underset{i}{0}, \dots, c_k, v) = 0$ & for every $i \geq 0$, \\
$\gamma_2 (c_1, \dots, \underset{i}{0}, \dots, c_k, v) = 0$ & for every $i \geq 1$. \\ 
\end{tabular}
\end{center}
It follows that $\widehat C_{\mathrm{def,lin}}(V_1)$ reduces to
\begin{equation} \label{eq:2-vect_comp}
0 \longrightarrow \mathrm{Hom}(V_0, C)[1] \overset{\delta_0}{\longrightarrow} \operatorname{End} C \oplus \operatorname{End} V_0 \overset{\delta_1}{\longrightarrow} \mathrm{Hom}(C, V_0)[-1] \longrightarrow 0
\end{equation}
with
\begin{align*}
\delta_0 \gamma & = (\gamma \circ \partial, \partial \circ \gamma), \\
\delta_1 (\gamma_1, \gamma_2) & = \gamma_2 \circ \partial - \partial \circ \gamma_1. \\
\end{align*}
So the linear deformation cohomology of $V_1$ is:
\begin{equation} \label{eq:2-vect_coh}
\begin{aligned}
H^{-1}_{\mathrm{def}, \mathrm{lin}} (V_1) & = \operatorname{Hom} (\operatorname{coker} \partial, \ker \partial),\\
H^0_{\mathrm{def}, \mathrm{lin}} (V_1) & = \operatorname{End} (\operatorname{coker} \partial) \oplus \operatorname{End} (\ker \partial), \\
H^1_{\mathrm{def}, \mathrm{lin}} (V_1) & = \operatorname{Hom}(\ker \partial, \operatorname{coker} \partial).
\end{aligned}
\end{equation}

Finally, notice that the VB-algebroid of $V_1$ is an \emph{LA-vector space} \cite{lapastina:def} of the form $(V_1 \Rightarrow V_0; 0 \Rightarrow *)$. In \cite{lapastina:def} we showed that the linear deformation complex of $(V_1 \Rightarrow V_0; 0 \Rightarrow *)$ is again (\ref{eq:2-vect_comp}), and it is easy to show, for example in coordinates, that the van Est map
\begin{equation*}
\mathrm{VE}: \widehat C_{\mathrm{def,lin}}(V_1 \rightrightarrows V_0) \to C_{\mathrm{def,lin}}(V_1 \Rightarrow V_0)
\end{equation*}
 is simply the identity.

\subsection{Tangent VB-groupoid}\label{sec:tangent-VB}

Let $\mathcal G \rightrightarrows M$ be a Lie groupoid. We want to relate the linear deformation cohomology of $T \mathcal G$ with the deformation cohomology of $\mathcal G$. First recall that there is a projection
\begin{equation}\label{eq:proj_TG}
\mathrm{pr}: C_{\mathrm{def,lin}}(T \mathcal G) \to C_{\mathrm{def}}(\mathcal G).
\end{equation}
Now define an inclusion
\[
\iota: C_{\mathrm{def}}(\mathcal G) \to C_{\mathrm{def,lin}}(T \mathcal G), \quad  c \mapsto \iota_c := J \circ Tc
\]
where $J: TT \mathcal G \to TT \mathcal G$ is the \emph{flip} of the DVB $TT \mathcal G$ (see, e.g., \cite{mackenzie} for definition and basic properties). For later use, we recall that, if $\pi: T\mathcal G \to \mathcal G$ is the canonical projection, the flip is an isomorphism of the DVB $TT \mathcal G$ that inverts the two projections $T \pi: TT \mathcal G \to T \mathcal G$ and $\pi_{T \mathcal G}: TT \mathcal G \to T \mathcal G$:
\begin{equation}\label{eq:J}
J \circ T \pi = \pi_{T\mathcal G} \circ J.
\end{equation}

The inclusion $\iota$ is well defined, i.e., if $c \in C^k_{\mathrm{def}}(\mathcal G)$, then $\iota_c \in C^k_{\mathrm{def,lin}}(T \mathcal G)$. To see this we show that properties (1) and (2) of Definition \ref{def:def_complex} hold. So, let $t \mapsto (g_0 (t), \dots, g_k (t))$ be a curve in $\mathcal G^{(k+1)}$ defined around 0. Then $(\pi \circ c)(g_0 (t), \dots, g_k (t)) = g_1 (t)$. Differentiating at $t = 0$ we get
\[
(T \pi \circ Tc)(\dot{g_0}(0), \dots, \dot{g_k}(0)) = \dot{g_0}(0).
\]
Applying $J$ and remembering equation (\ref{eq:J}) we get
\[
\pi_{T\mathcal G} (\iota_c (\dot{g_0}(0), \dots, \dot{g_k}(0))) = \dot{g_0}(0),
\]
i.e. $\iota_c (\dot{g_0}(0), \dots, \dot{g_k}(0)) \in T_{\dot{g_1}(0)} T \mathcal G$ as desired.

Now notice that property (2) of Definition \ref{def:def_complex} can be expressed in the following way: if $p: \mathcal G^{(k+1)} \to \mathcal G^{(k)}$ is the map that forgets the first arrow, then $T\mathsf s \circ c$ descends to a map $\mathcal G^{(k)} \to TM$:
\[
\begin{array}{c}
\xymatrix@C=45pt{
& \mathcal{G}^{(k+1)} \ar[r]^-{T \mathsf s \circ c} \ar[d]_-p & TM \\
& \mathcal G ^{(k)} \ar[ur]}
\end{array}.
\]
Differentiating, we obtain that also $TT \mathsf s \circ Tc: T \mathcal G^{(k+1)} \to TTM$ descends to a map $T \mathcal G^{(k)} \to TTM$. Applying $J$ and remembering that it commutes with $TT \mathsf s$, we obtain the same statement for $TT \mathsf s \circ \iota_c$, as desired. Moreover, $\iota_c$ is obviously linear.

A direct computation shows that $\iota$ is a cochain map. Finally, we observe that the following diagram commutes:
\[
\begin{array}{c}
\xymatrix{
 T \mathcal{G}^{(k+1)} \ar[r]^-{Tc} \ar[d] & TT \mathcal G \ar[d]^-{\pi_{T\mathcal G}} \ar[r]^-{J} & TT \mathcal G \ar[d]^-{T \pi} \\
 \mathcal G ^{(k+1)} \ar[r]^-{c} & T \mathcal G \ar@{=}[r] & T \mathcal G 
}
\end{array}.
\]
This shows that $\iota$ inverts the projection (\ref{eq:proj_TG}). It follows that
\[
C_{\mathrm{def,lin}}(T \mathcal G) \cong C_{\mathrm{def}}(\mathcal G) \oplus \ker (\mathrm{pr})
\]
as cochain complexes, hence
\[
H_{\mathrm{def,lin}}(T\mathcal G) = H_{\mathrm{def,lin}}(T^* \mathcal G) = H_{\mathrm{def}}(\mathcal G) \oplus H (\ker(\mathrm{pr})).
\]

\subsection{Representations of foliation groupoids}\label{sec:fol_grpd}

A \emph{foliation groupoid} is a Lie groupoid whose anchor map is injective. This condition ensures that the connected components of the orbits of the groupoid are the leaves of a regular foliation of the base manifold, whence the name. On the other hand, foliation groupoids encompass several classical groupoids associated to a foliated manifold, such as the holonomy and the monodromy groupoids.

By Example \ref{ex:act_grpd}, representations of a foliation groupoids $\mathcal G$ are equivalent to trivial core VB-groupoids over $\mathcal G$. Here we want to study the linear deformation cohomology of such VB-groupoids. First of all, let $\mathcal G \rightrightarrows M$ be a foliation groupoid, let $A \Rightarrow M$ be its Lie algebroid, $\rho: A \to TM$ the (injective) anchor map, $\nu = TM/ \operatorname{im} \rho$ the normal bundle and let $\pi: TM \to \nu$ be the projection. In this case, $\nu$ has constant rank and the normal representation is a plain representation of $\mathcal G$ on $\nu$. Moreover, we recall from \cite{crainic:def2} that the map
\[
p: C_{\mathrm{def}}(\mathcal G) \to C (\mathcal G, \nu), \quad c \mapsto \pi \circ s_c
\]
is a surjective quasi-isomorphism.

Consider a representation $E \to M$ of $\mathcal G$ and construct the associated trivial core VB-groupoid $(\mathcal G \ltimes E \rightrightarrows E; \mathcal G \rightrightarrows M)$. At the infinitesimal level, there is an induced representation of $A$ on $E$, i.e. an $A$-flat connection $\nabla: A \to DE$, and $\nabla$ is injective because $\rho$ is so. Therefore, the cokernel $\tilde \nu = DE/ \operatorname{im} \nabla$ is a vector bundle over $M$. Denote by
\[
\tilde \pi: DE \to \tilde \nu, \quad \delta \mapsto \bar \delta
\]
the projection.

We want to show that, in this situation, $\mathcal G$ acts on $\tilde \nu$. 
To see this, recall that the group $\mathrm{Bis}(\mathcal G)$ of bisections of $\mathcal G$ acts on $\Gamma(E)$ via
\begin{equation} \label{eq:bis_act_1}
(\beta \star \varepsilon)_x = \beta_{(\mathsf t \circ \beta)^{-1} (x)} \cdot \varepsilon_{(\mathsf t \circ \beta)^{-1} (x)},
\end{equation}
so it also acts on derivations of $E$ by
\begin{equation} \label{eq:bis_act_2}
(\beta \bullet \Delta)(\varepsilon) = \beta \star (\Delta( \beta^{-1} \star \varepsilon)).
\end{equation}
If $\beta$ is a local bisection around $x$ and $\beta(x) = g: x \to y$, these formulas still make sense: Equation (\ref{eq:bis_act_1}) shows that the action of $\beta$ takes local sections around $y$ to local sections around $x$, Equation (\ref{eq:bis_act_2}) shows that $\beta$ acts on derivations locally defined around $x$ (to give a derivation locally defined around $y$).

Now, let $g: x \to y$ be an arrow in $\mathcal G$ and $\delta \in D_x E$. Choose a local bisection $\beta$ of $\mathcal G$ passing through $g$ and a derivation $\Delta \in \mathfrak D (E)$ such that $\Delta_x = \delta$. Our action is then defined by
\[
g \cdot \bar \delta = \overline{\beta \bullet \Delta}|_y.
\]
A routine computation shows that the definition does not depend on the choice of $\Delta$. Let us prove that this definition is also independent of the choice of $\beta$. This is equivalent to prove that, if $\beta \in \mathrm{Bis}(\mathcal G)$, $\beta_z = \mathsf 1_z$ for some $z \in M$, then there exists $\alpha \in \Gamma(A)$ such that 
\[
(\Delta - \beta \bullet \Delta - \nabla_\alpha)_z = 0.
\]
Consider the vector field $\sigma_{\Delta - \beta \bullet \Delta}$. We have
\begin{equation}\label{eq:xi}
\xi := (\sigma_{\Delta - \beta \bullet \Delta})_z = (\sigma_\Delta - (\mathsf t \circ \beta)_* (\sigma_\Delta))_z = \sigma_{\Delta_z} - T(\mathsf t \circ \beta)(\sigma_{\Delta_z}).
\end{equation}
But $\mathsf t \circ \beta$ preserves the orbits of $\mathcal G$, hence it maps a sufficiently small neighborhood of $z$ in the leaf $\mathcal L_z$ of $\operatorname{im} \rho$ through $z$ to itself. It then follows from (\ref{eq:xi}) that $\xi$ kills all the functions that are constant along $\mathcal L_z$, hence it belongs to $\operatorname{im} \rho$.

Now, let $\alpha \in \Gamma(A)$ be any section such that $\rho (\alpha_z) = \xi$ and put $D := \Delta - \beta \bullet \Delta - \nabla_\alpha$. By construction, $\sigma_{D_z} = 0$, so it suffices to show that $D_z$ vanishes on $\nabla$-flat sections. If $\varepsilon$ is such a section, then
\begin{equation}\label{eq:D_z}
D_z \varepsilon = (\Delta - \beta \bullet \Delta)_z \varepsilon = \Delta_z (\varepsilon - \beta^{-1} \star \varepsilon).
\end{equation}
But the hypothesis $\nabla \varepsilon = 0$ implies that $\varepsilon$ is invariant under the action of $\mathrm{Bis} (\mathcal G)$, at least locally around $z$, and the claim follows from (\ref{eq:D_z}).

Notice that the symbol map $\sigma: DE \to TM$ descends to a $\mathcal G$-equivariant map $\tilde \nu \to \nu$, and the fact that $\operatorname  {End} E \cap \operatorname{im} \nabla = 0$ implies that its kernel is again $\operatorname{End} E$. Hence we have a short exact sequence of vector bundles with $\mathcal G$-action:
\[
0 \longrightarrow \operatorname{End} E \longrightarrow \tilde \nu \longrightarrow \nu \longrightarrow 0.
\]
In turn, this induces a short exact sequence of DG-modules:
\begin{equation} \label{eq:coch_1}
0 \longrightarrow C(\mathcal G, \operatorname{End} E) \longrightarrow C(\mathcal G, \tilde \nu) \longrightarrow C(\mathcal G, \nu) \longrightarrow 0.
\end{equation}
But $\mathcal G \ltimes E$ is also a trivial-core VB-groupoid, so we also have the sequence (\ref{eq:coch_2}):
\[
0 \longrightarrow C(\mathcal G, \operatorname{End} E) \longrightarrow C_{\mathrm{def,lin}}(\mathcal G \ltimes E) \longrightarrow C_{\mathrm{def}}(\mathcal G) \longrightarrow 0.
\] 
We are looking for a map relating the two sequences. If $\tilde c \in C^k_{\mathrm{def,lin}}(\mathcal G \ltimes E)$ and $\tilde c_2: \mathcal G^{(k)} \to DE$ is the map (\ref{eq:c_2}), we can simply define:
\[
\tilde p: C_{\mathrm{def,lin}}(\mathcal G \ltimes E) \to C(\mathcal G, \tilde \nu), \quad \tilde c \mapsto \tilde \pi \circ \tilde c_2.
\]
This is a cochain map and we obtain the following commutative diagram:
\[
\begin{array}{c}
\xymatrix@C=15pt@R=20pt{
	0 \ar[r] & C(\mathcal G, \operatorname{End} E) \ar[r] \ar@{=}[d] & C_{\mathrm{def,lin}}(\mathcal G \ltimes E) \ar[r] \ar[d]^{\tilde p} & C_{\mathrm{def}}(\mathcal G) \ar[r]  \ar[d]^{p}  & 0 \\
	0 \ar[r] & C(\mathcal G, \operatorname{End} E) \ar[r] & C(\mathcal G, \tilde \nu) \ar[r] & C(\mathcal G, \nu) \ar[r] & 0
}
\end{array}.
\]
The rows are short exact sequences of DG-modules and the vertical arrows are DG-module surjections; additionally, $p$ is a quasi-isomorphism. Hence, it immediately follows from the Snake Lemma and the Five Lemma that $\tilde p$ is a quasi-isomorphism as well. We have thus proved the following
\begin{proposition} There is a canonical isomorphism between the linear deformation cohomology of the VB-groupoid $(\mathcal G \ltimes E \rightrightarrows E; \mathcal G \rightrightarrows M)$ and the leafwise cohomology with coefficients in $\tilde \nu$:
	\[
	H_{\mathrm{def}, \mathrm{lin}} (\mathcal G \ltimes E) = H (\mathcal G, \tilde \nu).
	\]
\end{proposition}

\subsection{Lie group actions on vector bundles}\label{sec:Lie_vect}

Let $G$ be a Lie group with Lie algebra $\mathfrak g$. Assume that $G$ acts on a vector bundle $E \to M$ by vector bundle automorphisms. Then $G$ acts also on $M$ and $(G \ltimes E \rightrightarrows E; G \ltimes M \rightrightarrows M)$ is a trivial-core VB-groupoid. We want to discuss its linear deformation cohomology.

Of course the action of $G$ on $M$ induces an infinitesimal action of $\mathfrak g$ on $M$, and the Lie algebroid of $G \ltimes M$ is the action algebroid $\mathfrak g \ltimes M$. We recall from \cite{crainic:def2} that $G \ltimes M$ acts naturally on $\mathfrak g \ltimes M$, by extending the adjoint action of $G$ on $\mathfrak g$, and on $TM$ by differentiating the action of $G$ on $M$. Moreover, there is a short exact sequence of complexes:
\begin{equation}\label{eq:act_seq}
0  \longrightarrow C (G \ltimes M, TM) \longrightarrow C_{\mathrm{def}}(G \ltimes M)  \overset{p}{\longrightarrow} C (G \ltimes M, \mathfrak g \ltimes M)[1] \longrightarrow 0.
\end{equation}
The projection $p$ is defined as follows. Take a $k$-cochain $c: (G \ltimes M)^{(k+1)} \to T(G \ltimes M)$ in the deformation complex of $G \ltimes M$ and let $(h_0, \dots, h_k) \in (G \ltimes M)^{(k+1)}$, $h_0 = (g, x)$. Then
\[
c(h_0, \dots, h_k) \in T_{(g_0,x_0)}(G \ltimes M) \cong T_{g_0} G \times T_{x_0} M \cong \mathfrak g \times T_{x_0} M
\]
via right translations, and we compose with the projection $TM \to M$ to get an element in $\mathfrak g \ltimes M$. The kernel of $p$ is given by $TM$-valued cochains. Since the $TM$-component is the projection by the source map, it does not depend on the first component, and we conclude that the kernel is $C(G \ltimes M, TM)$.

We want to construct a similar sequence for $C_{\mathrm{def,lin}}(G \ltimes E)$ taking into account the linear nature of the action.  First of all, there is an obvious induced action of $G$ on $DE$ and the symbol map $\sigma: DE \to TM$ is $G$-equivariant. Hence there is a short exact sequence of cochain complexes:
\[
0 \longrightarrow C(G \ltimes M, \operatorname{End} E) \longrightarrow C(G \ltimes M, DE) \longrightarrow C(G \ltimes M, TM) \longrightarrow 0.
\]
 In this case, the sequence (\ref{eq:coch_2}) reads
\begin{equation}\label{eq:Lie_vect_pr_0}
0 \longrightarrow C(G \ltimes M, \operatorname{End} E) \longrightarrow C_{\mathrm{def,lin}}(G \ltimes E) \overset{\Pi}{\longrightarrow}  C_{\mathrm{def}} (G \ltimes M) \longrightarrow 0
\end{equation}
and, composing $\Pi$ with $p$, we get a cochain map
\begin{equation}\label{eq:Lie_vect_pr}
C_{\mathrm{def,lin}}(G \ltimes E) \longrightarrow C (G \ltimes M, \mathfrak g \ltimes M)[1] \longrightarrow 0.
\end{equation}
Applying the isomorphism (\ref{eq:nerve_act}), we get $(G \ltimes E)^{(k)} \cong G^k \times E$, and similarly $(G \ltimes M)^{(k)} \cong G^k \times M$. Then if $\tilde c \in C^k_{\mathrm{def,lin}}(G \ltimes E)$ is a linear cochain, the diagram (\ref{diag:c_tilde_2}) takes the following form:
\[
\begin{array}{c}
\xymatrix@C=0pt@R=7pt{
& G^{k+1} \times E \ar[dl] \ar[rr]^{\tilde c} \ar[dd]|!{[dl];[dr]}{\hole} & & TG \times TE \ar[dl] \ar[dd] \\
G^k \times E \ar[rr]^(.65){s_{\tilde c}} \ar[dd] & & TE \ar[dd] \\
& G^{k+1} \ar[dl] \ar[rr]^(.4)c |!{[ur];[dr]}{\hole} & & TG \times TM  \ar[dl] \\
G^k \times M \ar[rr]^(.65){s_c} & & TM  }
\end{array}.
\]
Moreover, $\tilde c$ is in the kernel of (\ref{eq:Lie_vect_pr}) if and only if its $TG$-component is $0$. In this case, it is clear that $\tilde c$ is determined by $s_{\tilde c}$, that is in turn equivalent to the map $\tilde c_2: G^k \times M \to DE$ defined by (\ref{eq:c_2}). Therefore, the kernel of (\ref{eq:Lie_vect_pr}) is $C(G \ltimes M, DE)$.

Summarizing there is an exact diagram of cochain complexes
\[
\begin{array}{c}
\xymatrix@C=15pt@R=15pt{
& 0 \ar[d] & 0 \ar[d] & & \\
0 \ar[r] & C (G \ltimes M, \operatorname{End} E) \ar[d] \ar@{=}[r] & C (G \ltimes M, \operatorname{End} E) \ar[r] \ar[d] & 0 \ar[d] \\
0 \ar[r] & C (G \ltimes M, DE) \ar[r] \ar[d] & C_{\mathrm{def}, \mathrm{lin}} (G \ltimes E) \ar[r] \ar[d] & C (G \ltimes M, \mathfrak g \ltimes M)[1] \ar[r] \ar@{=}[d] & 0 \\
0 \ar[r] & C (G \ltimes M, TM) \ar[r] \ar[d] & C_{\mathrm{def}}(G \ltimes M) \ar[r] \ar[d] & C (G \ltimes M, \mathfrak g \ltimes M)[1] \ar[r] \ar[d] & 0 \\
 & 0 & 0 & 0 &
}
\end{array}.
\]
This proves the following
\begin{proposition}
Let $G$ be a Lie group acting on a vector bundle $E \to M$ by vector bundle automorphisms. The linear deformation cohomology of the VB-groupoid $(G \ltimes E \rightrightarrows E, G \ltimes M \rightrightarrows M)$ fits in the exact diagram:
\[
\begin{array}{c}
\xymatrix@C=15pt@R=15pt{
& \vdots \ar[d] &  \vdots \ar[d] &  \vdots \ar[d] & \\
\cdots \ar[r] & H^{k} (G \ltimes M, \operatorname{End} E) \ar@{=}[r] \ar[d]& H^{k}(G \ltimes M, \operatorname{End} E) \ar[r] \ar[d]& 0 \ar[r] \ar[d] & \cdots \\
\cdots \ar[r] & H^{k} (G \ltimes M, DE) \ar[r] \ar[d]& H^k_{\mathrm{def}, \mathrm{lin}}(G \ltimes E) \ar[r] \ar[d]& H^{k+1} (G \ltimes M, \mathfrak g \ltimes M) \ar[r] \ar@{=}[d] & \cdots \\
\cdots \ar[r] & H^{k} (G \ltimes M, TM) \ar[r] \ar[d] & H^k_{\mathrm{def}}(G \ltimes M) \ar[r] \ar[d]& H^{k+1} (G \ltimes M, \mathfrak g \ltimes M) \ar[r] \ar[d] & \cdots \\
\cdots \ar[r] & H^{k+1} (G \ltimes M, \operatorname{End} E) \ar@{=}[r] \ar[d]& H^{k+1}(G \ltimes M, \operatorname{End} E) \ar[r] \ar[d]& 0 \ar[r] \ar[d] & \cdots \\
& \vdots  &  \vdots  &  \vdots  &
}
\end{array}.
\]
\end{proposition}

\appendix

\section{} \label{sec:app}

In this section we prove a technical result that is needed in Subsection \ref{sec:lin_def}. By a \emph{fibration of DVBs} we mean a DVB morphism 
\[
\xymatrix@C=20pt@R=15pt{
& W_1 \ar[dl] \ar[rr] \ar[dd]|!{[dl];[dr]}{\hole} & & W_2 \ar[dl] \ar[dd] \\
E_1 \ar[rr] \ar[dd] & & E_2 \ar[dd] \\
& A_1 \ar[dl] \ar[rr]|!{[ur];[dr]}{\hole} & & A_2 \ar[dl] \\
M_1 \ar[rr] & & M_2 }
\]
where all the horizontal maps are surjective submersions.

\begin{lemma}\label{lemma:section} Consider a fibration of DVBs:
\[
\xymatrix@C=20pt@R=15pt{
& W_1 \ar[dl] \ar[rr]^(.6){\phi} \ar[dd]|!{[dl];[dr]}{\hole} & & W_2 \ar[dl] \ar[dd] \\
E_1 \ar[rr]^(.65){\phi_E} \ar[dd] & & E_2 \ar@/_0.7pc/[ur]|(.5){\tilde \alpha_2} \\
& A_1 \ar[dl] \ar[rr]|!{[ur];[dr]}{\hole} & & A_2 \ar[dl] \\
M_1 \ar[rr] \ar@/_0.7pc/[ur]|(.6){\alpha_1} & & M_2 \ar[uu] \ar@/_0.7pc/[ur]|(.6){\alpha_2} }
\]
let $\tilde \alpha_2$ be a linear section of $W_2 \to E_2$, and let $\alpha_1$ be a section of $A_1$ that projecting on the same section $\alpha_2$ of $A_2$. Then there exists a linear section of $W_1 \to E_1$ that projects simultaneously on $\tilde \alpha_2$ and $\alpha_1$. \end{lemma}

\proof We can decompose the morphism $\phi$ in the following way:
\[
\xymatrix@C=15pt@R=15pt{
& W_1 \ar[dl] \ar[rr]^-{\widehat \phi} \ar[dd]|!{[dl];[dr]}{\hole} & & E_1 \times_{E_2} W_2 \times_{A_2} A_1 \ar[dl] \ar[dd]|!{[dl];[dr]}{\hole} \ar[rr] & & W_2 \ar[dd] \ar[dl] \\
E_1 \ar@{=}[rr] \ar@/_0.7pc/[ur]|(.6){\tilde \alpha_1} \ar[dd]^-{\pi_{E_1}} & & E_1 \ar@/_0.7pc/[ur]|(.55){\widehat{\alpha}} \ar[rr]^(.6){\phi_E} \ar[dd] & & E_2 \ar[dd] \ar@/_0.7pc/[ur]|-{\tilde \alpha_2} \\
& A_1 \ar[dl] \ar@{=}[rr]|!{[ur];[dr]}{\hole} & & A_1 \ar[dl] \ar[rr]|!{[ur];[dr]}{\hole} & & A_2 \ar[dl] \\
M_1 \ar@{=}[rr] & & M_1 \ar@/_0.7pc/[ur]|(.6){\alpha_1} \ar[rr] & & M_2 \ar@/_0.7pc/[ur]|-{\alpha_2}
}
\]
Define
\[
\widehat \alpha: E_1 \to E_1 \times_{E_2} W_2 \times_{A_2} A_1, \quad \widehat \alpha_{e} = (e, \tilde \alpha_2 |_{\phi_E(e)}, \alpha_1 |_{\pi_{E_1}(e)}).
\]
It is clear that $\widehat \alpha$ is a well-defined linear section of $E_1 \times_{E_2} W_2 \times_{A_2} A_1 \to E_1$, that projects simultaneously on $\tilde \alpha_2$ and $\alpha_1$, so the problem is reduced to find a splitting of $\widehat \phi$. This can be done first locally, and then globally via the choice of a partition of unity on $M_1$. \endproof

\bigskip

\noindent \textbf{Acknowledgments.} The authors are members of GNSAGA of INdAM. This work has been partially carried out while the first author was visiting the Department of Mathematics and Statistics of S\~{a}o Paulo University, with the financial support of INdAM. He would like to thank Cristian Ortiz and Ivan Struchiner for many fruitful discussions and the all department for the kind hospitality.

\end{document}